\documentclass[11pt]{amsart}
\usepackage{latexsym,enumerate}
\usepackage{amsmath,amsthm,amsopn,amstext,amscd,amsfonts,amssymb}
\usepackage{color}
\numberwithin{equation}{section}

\newtheorem{theorem}{Theorem}
\newtheorem{lemma}{Lemma}
\newtheorem{corollary}{Corollary}
\newtheorem{proposition}{Proposition}
\newtheorem{remark}{Remark}
\newtheorem{definition}{Definition}

\numberwithin{theorem}{section}
\numberwithin{corollary}{section}
\numberwithin{lemma}{section}
\numberwithin{definition}{section}
\numberwithin{proposition}{section}
\numberwithin{remark}{section}

\newcommand{\RR}{\mathbb R^N}

\newcommand{\medint}{-\kern  -,375cm\int}

\begin{document}
\title[ Weighted isoperimetric inequalities in cones and applications ]
{  Weighted isoperimetric inequalities in cones \\ and applications }
\author{ F. Brock$^1$ - F. Chiacchio$^2$ - A. Mercaldo$^2$}
\thanks{}
\date{}

\begin{abstract}
\noindent This paper deals with weighted isoperimetric 
inequalities relative to cones of $\mathbb{R}^{N}$.
We study the structure of measures that admit as isoperimetric 
sets the intersection of a cone
with balls centered at the vertex of the cone.
For instance, in case that the cone is the half-space
$\mathbb{R}_{+}^{N}=\left\{ x \in \mathbb{R}^{N} : x_{N}>0 \right\}$
and the measure is factorized, we prove that this phenomenon occurs if and only if
the measure has the form $d\mu=ax_{N}^{k}\exp 
\left( c\left\vert x\right\vert ^{2}\right)dx $, for some $a>0$, $k,c\geq 0$.
Our results are then used to obtain isoperimetric estimates for
Neumann eigenvalues of a weighted Laplace-Beltrami operator on the
sphere, sharp Hardy-type inequalities for functions defined in a
quarter space and, finally, via symmetrization arguments, a
comparison result for a class of degenerate PDE's.
\bigskip

\noindent \textsl{Key words:} \textrm{relative isoperimetric inequalities,
Neumann eigenvalues, weighted Laplace-Beltrami operator, Hardy inequalities, 
degenerate elliptic equations. \medskip }

\end{abstract}

\maketitle

\setcounter{footnote}{1} \footnotetext{
American University of Beirut, Department of Mathematics, Beirut, Lebanon,
P.O. Box: 11-0236, e-mail: fb13@aub.edu.lb, 
partially supported by FONDECYT project 1050412}

\setcounter{footnote}{2} \footnotetext{
Dipartimento di Matematica e Applicazioni \textquotedblleft R.
Caccioppoli\textquotedblright , Universit\`{a} degli Studi di Napoli
``Fe\-derico II", Complesso Monte S. Angelo, via Cintia, 80126 Napoli,
Italy, e-mails: francesco.chiacchio@unina.it, mercaldo@unina.it}

\section{Introduction}

This paper deals with weighted relative isoperimetric inequalities in cones
of $\mathbb{R}^{N}$. Let $\omega $ be an open subset of  $\mathbb{S}^{N-1} $,
the unit sphere of $\mathbb{R}^{N}$, and
$\Omega $ the cone
\begin{equation}
\label{introcone}
\Omega =\left\{ x\in \mathbb{R}^{N}:\frac{x}{\left\vert x\right\vert }\in
\omega , x\not=0  \right\} .
\end{equation}
We consider measures of the type $d\nu =\phi (x)dx$ on $\Omega $,
where $\phi $ is a positive Borel measurable function defined in $\Omega $.
For any measurable set $M\subset \Omega $, we define the $\nu $
-measure of $M$
\begin{equation}
\label{intronu}
\nu (M)=\int_{M}d\nu =\int_{M}\phi(x)dx
\end{equation}
and the $\nu $-perimeter of $M$ relative to $\Omega$
\begin{equation*}
P_{\nu }(M,\Omega )=\sup \left\{  \int_{M}\mbox{div }(\mathbf{v}(x)\phi(x))\,dx:\,
\mathbf{v}\in C_0 ^1(\Omega ,\mathbb{R}^{N}),\,|\mathbf{v}|\leq 1\    \right\}\, .
\end{equation*}
We also write $P_{\nu }(M,\mathbb{R}^{N})=P_{\nu }(M)$. Note that if $M$
is a smooth set, then
\begin{equation*}
P_{\nu }(M,\Omega )=\int_{\partial M\cap \Omega }\phi(x)\, d{\mathcal H}_{N-1}(x).
\end{equation*}
The isoperimetric problem reads as
\begin{equation}\label{Imu}
I_{\nu }(m)=\inf \{P_{\nu }(M,\Omega ):\, M\subset \Omega ,\,\nu
(M)=m\},\quad m>0.
\end{equation}
One says that $M $ is an isoperimetric set if $\nu (M)=m$ and $I_{\nu }(m)=
P_{\nu }(M,\Omega )$.

We give necessary conditions on the function $\phi$ for having
$B_{R}\cap \Omega $ as an isoperimetric set, in Section 2.
Here and
throughout the paper, $B_R$ and $B_R(x)$ denote the ball of radius $R$
centered at zero and at $x$, respectively.
In Theorem 2.1 we prove
that if $B_{R}\cap \Omega $ is an isoperimetric set for every $R>0$,  then
\begin{equation*}
\phi=A(r)B(\Theta ),
\end{equation*}
where $r=\left\vert x\right\vert $ and $\Theta =\frac{x}{\left\vert
x\right\vert }.$
\\
As an application of Theorem 2.1, we prove a sharp Hardy-type inequality for
functions defined in $Q =\left\{ x_{1}>0, x_{N}>0 \right\}$
involving a power-type weight, (see Theorem 2.6). 

We are able to give an explicit expression of the density 
$\phi $ in some special cases.
For instance, 
when $\Omega $ is the   
half space
\begin{equation}
\Omega =\mathbb{R}_{+}^{N}=\left\{ x=(x_{1},...,x_{N})\in \mathbb{R}
^{N}:x_{N}>0\right\} ,  \label{R+} 
\end{equation}
if $\phi$ is a smooth function with a factorized structure,
\begin{equation}
\phi(x)=\displaystyle\prod\limits_{i=1}^{N}\phi_{i}(x),  
\label{intro_f_fact}
\end{equation}
and if $B_{R}\cap \mathbb{R}_{+}^{N} $ is an isoperimetric set,
then
\begin{equation}
\phi(x)=ax_{N}^{k}\exp( c\left\vert x\right\vert ^{2}) ,
\label{intro_f}
\end{equation}
for some numbers $a>0,$ $k\geq 0$ and $c\geq 0$, (see Theorem 2.8).

Section 3 is dedicated to the case $\Omega = \mathbb{R}^N $, 
and to the proof of the following Theorem,
which is the main result of our paper.
%\\[0.2cm]

\begin{theorem}
\label{Theorem 1.1.}
{\sl Let  $\mu$ be the measure defined by
\begin{equation}
d\mu =x_{N}^{k}\exp ( c\left\vert x\right\vert ^{2}) dx,\text{ \ }
x\in \mathbb{R}_{+}^{N},  
\label{intro_dmu}
\end{equation}
with $k,c\geq 0$, and let $M$ be a measurable subset of 
$\mathbb{R}^N _+ $ with finite 
$\mu$-measure. Then
$$
P_{\mu }(M )\geq P_{\mu }(M ^{\bigstar })\, ,
$$
where $M ^{\bigstar }=B_{r^{\bigstar }}\cap \mathbb{R}_{+}^{N}$, with
$r^{\bigstar }$ such that $\mu (M )=\mu (M ^{\bigstar })$.}
\end{theorem}
%\\[0.2cm]

The proof of Theorem 1.1 requires some technical effort which is due to the degeneracy of the measure on the hyperplane $\{ x_N =0 \} $.

Note that
Theorem 1.1 is imbedded in a wide bibliography related to the
isoperimetric problems for \textquotedblleft manifolds with
density\textquotedblright\ (see, for instance, \cite{Bo}, \cite{BMP},
\cite{CMV}, \cite{CJQW}, \cite{DDNT}, \cite{Mo}, \cite{Mo2}, \cite{RCBM}, 
\cite{S}). Further references will be given in Section 2.

It was shown in \cite{MS} that the isoperimetric set
for measures of the type $y^{k}dxdy$, with $k\geq 0$ and $(x,y)\in \mathbb{R}_{+}^{2}$, is
$B_{R}\cap \mathbb{R}_{+}^{2}$.  In \cite{Borell} C. Borell  
 proved that balls centered at the origin are isoperimetric sets for measures of the type $\exp (
c\left\vert x\right\vert ^{2}) dx$ in $\mathbb{R}^{N}$ with $c\geq 0$ 
(see also \cite{BMP} and \cite{RCBM} for this and related results).

In Section 4 we consider
 degenerate
elliptic problems of the type
\begin{equation}
\left\{
\begin{array}{ccc}
-\text{div}(A(x)\nabla u)= & x_{N}^{k}\exp (
c\left\vert x\right\vert
^{2}) f(x) & \text{in }D \\
&  &  \\
u=0\text{ \ \ \ \ \ \ \ \ \ \ \ } &  & \text{ \ on }\Gamma _{+,}
\end{array}
\right.  
\label{P}
\end{equation}
where $D $ is a bounded open set in $\mathbb{R}_{+}^{N},$ whose
boundary is decomposed into a part $\Gamma _{0}$, lying on the
hyperplane $\left\{ x_{N}=0\right\} $ and a part $\Gamma _+$
contained in $\mathbb{R}_{+}^{N}$. (For precise definitions, see Section 4).
Assume that $c,k\geq 0,$ $A(x)=(a_{ij}(x))_{ij}$
is an $N\times N$ symmetric matrix with measurable
coefficients satisfying
\begin{equation}
x_{N}^{k}\exp ( c\left\vert x\right\vert ^{2}) \left\vert \zeta
\right\vert ^{2}\leq a_{ij}(x)\zeta _{i}\zeta _{j}\leq \Lambda x_{N}^{k}\exp
( c\left\vert x\right\vert ^{2})
\left\vert \zeta \right\vert
^{2},\text{ \ \ }\Lambda \geq 1,  
\label{ell}
\end{equation}
for almost every $x\in D $ \ and for all $\zeta \in \mathbb{R}^{N}$. Assume also
that $f$ belongs to the weighted Lebesgue space $L^{2}(D
,d\mu )$ where $d\mu $\ is the measure defined in (\ref{intro_dmu}).

The type of degeneracy in (\ref{ell}) occurs, for $k\in \mathbb{N}$, when
one looks for solutions to linear PDE's which are symmetric with respect to a group of $(k+1)$
variables (see, e.g., \cite{BCM}, \cite{MS}, \cite{T0}
and the references therein). The case of a non-integer $k$ has been the object of
investigation, for instance, in the generalized axially symmetric potential theory
(see, e.g., \cite{W} and the subsequent works of A. Weinstein).

We obtain optimal bounds for the solution to problem (\ref{P})
using a  symmetrization technique which is
due to G. Talenti (see \cite{T2} and also \cite{AFLT},
\cite{BBMP1}, \cite{BBMP3}, \cite{BCM}, \cite{MS}, \cite{RV}).

If $M$ is measurable set with finite $\mu$-measure, and if 
$f: M \rightarrow \mathbb{R} $ is a measurable function, the
 weighted rearrangement
 $f^{\bigstar }:M
^{\bigstar }\rightarrow \left[ 0,+\infty \right[ $ is uniquely defined by
the following condition
\begin{equation}\label{defrearr}
\left\{ x\in M ^{\bigstar }:f^{\bigstar }(x)>t\right\} =\left\{ x\in
M :\left\vert f(x)\right\vert >t\right\} ^{\bigstar }\text{ \ }\forall
t\geq 0.
\end{equation}
This means that the super  level sets of $f^{\bigstar }$ are half-balls centered
at the origin, having the same
$\mu- $measure of the corresponding super level sets of
$\left\vert f\right\vert $.

Let $C_{\mu }$ denote the $\mu $-measure of $B_{1}\cap \mathbb{R}_{+}^{N}$.
Using Theorem 1.1, we obtain the following comparison result.
%\\[0.2cm]
\begin{theorem}\label{Theorem 1.2.} {\sl 
Let $u$ be the weak solution to problem (\ref{P}), and  let $w$ be
the function
\begin{equation*}
w(x)=w^{\bigstar }(x)=\frac{1}{{C_{\mu }}}\int_{\left\vert x\right\vert
}^{r^{\bigstar }}\left( \int_{0}^{\rho }f^{\bigstar }(\sigma )\sigma
^{N-1+k}\exp \left( c\sigma ^{2}\right) d\sigma \right) \rho ^{-N+1-k}\exp
\left( -c\rho ^{2}\right) d\rho ,
\end{equation*}
which is the weak solution to the problem
\begin{equation}
\left\{
\begin{array}{ccc}
-\text{div}\left( x_{N}^{k}\exp \left( c\left\vert x\right\vert ^{2}\right)
\nabla w\right) =x_{N}^{k}\exp \left( c\left\vert x\right\vert ^{2}\right)
f^{\bigstar } & \text{in} & D ^{\bigstar } \\
&  &  \\
w=0\text{ \ \ \ \ \ \ \ \ \ \ \ \ \ \ \ \ \ \ \ \ \ \ \ \ \ \ \ \ \ \ \ \ \
\ \ \ \ \ \ \ \ \ \ \ \ \ } & \text{on} & \partial D ^{\bigstar }\cap
\mathbb{R}_{+}^{N}.
\end{array}
\right.  
\label{Symm_Probl_dv}
\end{equation}
Then
\begin{equation}
u^{\bigstar }(x)\leq w(x)\text{ a.e. in }D ^{\bigstar },
\label{Point_Est-dv}
\end{equation}
and
\begin{equation}
\int_{D }\left\vert \nabla u\right\vert ^{q}d\mu \leq \int_{D ^{\bigstar
}}\left\vert \nabla w\right\vert ^{q}d\mu ,\text{ for all }0<q\leq 2,
\label{Grad-Est-dv}
\end{equation}
}
\end{theorem}

\section{Weighted isoperimetric inequalities in a cone of $\mathbb{R}^{N}$}
In this section we study isoperimetric problems with respect to
measures, relative to cones in $\mathbb{R}^{N}$. Notice that such
problems have been investigated for instance in \cite{ACDLV}, \cite{Ba},
\cite{DHHT}, \cite{LP}, \cite{MR} and \cite{PT}. Our aim is to
characterize those measures for which an isoperimetric set is given
by the intersection of a cone with the ball having center at the
vertex of the cone.

We begin by fixing some notation that will be used throughout:
$\omega _{N}$ is the $N$-dimensional Lebesgue measure of the unit
ball in $\mathbb{R}^{N}$. For points $x\in \mathbb{R}^{N}-\{0\}$ we will
often use $N$-dimensional polar coordinates $(r,\Theta )$, where
$r=|x|$ and $\Theta =x|x|^{-1}\in \mathbb{S} ^{N-1}$. $\nabla
_{\Theta } $ denotes the gradient on $\mathbb{S} ^{N-1} $.
 By $\mathbb{S}_+ ^{N-1} $ we denote the half sphere,
$$
\mathbb{S}_+ ^ {N-1} = \mathbb{S} ^{N-1} \cap \mathbb{R} _+ ^{N} .
$$
Consider the isoperimetric problem \eqref{Imu} where $\Omega$ is the cone
defined in \eqref{introcone} and $\nu$ the measure given by \eqref{intronu}.

The first result of this section
says that, if the isoperimetric set of (\ref{Imu}) is $B_{R}\cap
\Omega $ for a suitable $R$, then the density of the measure $d\nu $
is a product of two functions $A$ and $B$ of the variables $r$ and
$\Theta $, respectively.

\noindent Note that it has been proven in \cite{LMPW} that a smooth density on $\RR$ is radial if and only if spheres about the origin are stationary for a given volume.
%\\[0.2cm]
\begin{theorem}\label{Theorem 2.1.} {\sl
Consider Problem (\ref{Imu}), with $\phi\in C^{1}(\Omega
)\cap C( \overline{\Omega })$, $\phi(x)>0$ for $x\in \Omega $. Suppose
that $I_{\nu }(m)=P_{\nu }(B_{R}\cap \Omega )$ whenever $m=\nu (B_{R}\cap
\Omega )$, for every $R>0$. Then
\begin{equation}
\phi =A(r)B(\Theta ),  
\label{decomp}
\end{equation}
where $A\in C^{1}((0,+\infty ))\cap C([0,+\infty ))$, $A(r)>0$ if $r>0$, and
$B\in C^{1}(\omega )$, $B(\Theta )>0$ for $\Theta \in \omega $. Moreover, if
$\phi \in C^{2}(\Omega )$, then
\begin{equation}
\lambda (B,\omega )\geq N-1+r^{2}\left[ \frac{(A^{\prime
} (r))^{2}}{ (A(r))^{2}}-\frac{A^{\prime \prime }(r)}{A(r)}\right] \quad
\forall r>0\,, \label{stability}
\end{equation}
where
\begin{equation}
\lambda (B,\omega ):=\inf \left\{ \frac{\int_{\omega }|\nabla _{\Theta
}u|^{2}B\,d\Theta }{\int_{\omega }u^{2}B\,d\Theta }:\,u\in C^{1}(\omega
),\,\int_{\omega }uB\,d\Theta =0,\,u\not=0\right\} \,.  \label{lambda1}
\end{equation}
}
\end{theorem}

\begin{remark}\label{Remark 2.2.}\rm
Observe that $\lambda (B,\omega )$ is the first
nontrivial eigenvalue of the Neumann problem
\begin{equation}
\left\{
\begin{array}{ccc}
-\nabla _{\Theta }\left( B\nabla _{\Theta }u\right) =\lambda Bu &\mbox{ in } &
\omega  
\\
&  &  
\\
\dfrac{\partial u}{\partial {\bf n} }=0 \text{ 
\ \ \ \ \ \ \ \ \ \ \
\ \ \ \ \ } & \mbox{ on } & \partial \omega
\end{array}
\right.   \notag
\end{equation}
where $u\in W^{1,2}(\omega )$, and ${\bf n} $ is the exterior unit
normal to
$ \partial \omega .$ 
\end{remark}

%\\[0.2cm]
\noindent{\sl Proof of Theorem 2.1 : }
Let $R>0$. For $\varepsilon \in \mathbb{ R}$ we define the following
measure-preserving perturbations $G_{\varepsilon }$ from $B_{R}\cap
\Omega $:
\begin{equation*}
G_{\varepsilon }:=\{(r,\Theta ):\,0<r<R+\varepsilon h(\Theta )+s(\varepsilon
),\,\Theta \in \omega \},\quad |\varepsilon |\leq \varepsilon _{0}
\end{equation*}
where $h\in C^{1}(\overline{\omega })$, and $s$ is to be chosen such that $s\in
C^{2}([-\varepsilon _{0},\varepsilon _{0}])$, $s(0)=0$, and $\nu
(G_{\varepsilon })=\nu (B_{R})$ for $|\varepsilon |\leq \varepsilon _{0}$.
Writing $\phi =\phi (r,\Theta )$, and
\begin{equation*}
R^{\varepsilon }:=R+\varepsilon h+s(\varepsilon ),
\end{equation*}
we have, for $|\varepsilon |\leq \varepsilon _{0}$,
\begin{equation}
\nu (G_{\varepsilon })=\int_{\omega }\int_{0}^{R^{\varepsilon }}r^{N-1}\phi
(r,\Theta )\,dr\,d\Theta =\nu (B_{R})
\label{measurepreserve}
\end{equation}
and
\begin{equation}
P_{\nu }(G_{\varepsilon },\Omega )=\int_{\omega }(R^{\varepsilon
})^{N-2}\phi (R^{\varepsilon }, \Theta )\sqrt{(R^{\varepsilon
})^{2}+|\nabla _{\Theta }R^{\varepsilon }|^{2}}\,d\Theta \geq P_{\nu
}(B_{R}\cap \Omega ,\Omega ). 
\label{minimize}
\end{equation}
Denote $s_{1}:=s^{\prime }(0)$ and
$ s_{2}:=s^{\prime \prime }(0)$.
Differentiating (\ref{measurepreserve}) gives
\begin{equation}
0=\int_{\omega }\phi (R,\Theta )(h(\Theta )+s_{1})\,d\Theta ,
\label{measurepreserve1}
\end{equation}
and
\begin{equation}
0=\int_{\omega }\left( (N-1)\phi (R,\Theta )+R\phi _{r}(R,\Theta )\right)
(h(\Theta )+s_{1})^{2}\,d\Theta +s_{2}R\int_{\omega }\phi (R,\Theta
)\,d\Theta .
\label{measurepreserve2}
\end{equation}
Using (\ref{minimize}) we get
\begin{equation}
\left\{
\begin{array}{c}
\displaystyle\frac{\partial }{\partial \varepsilon} P_{\nu
}(G_{\varepsilon },\Omega )
)\Big| _{\varepsilon =0}=0 \\
\\
\displaystyle\frac{ \partial ^{2}}{\partial \varepsilon ^2 } P_{\nu
}(G_{\varepsilon },\Omega ) \Big| _{\varepsilon =0}\geq 0.
\end{array}
\right.
\label{diffPmu_1}
\end{equation}
The first condition in (\ref{diffPmu_1}) gives
\begin{equation}
\int_{\omega }\left( (N-1)\phi (R,\Theta )+R\phi _{r}(R,\Theta )\right)
(h(\Theta )+s_{1})\,d\Theta =0.  
\label{diff1P}
\end{equation}
In other words, we have that $\int _{\omega } ( (N-1) \phi + R\phi _r )v \, d\theta =0 $ for all functions $v\in C^1
(\overline{\omega })$ satisfying $\int _{\omega } \phi v \, d\theta =0 $. Then the Fundamental
Lemma in the Calculus of Variations tells us that there is a
number $k(R)\in \mathbb{R}$ such that
\begin{equation}
\phi _{r}(R,\Theta )=k(R)\phi (R,\Theta )\quad \forall \Theta \in
\omega .
\label{lambdaR}
\end{equation}
Integrating this with respect to $R$ implies (\ref {decomp}). Hence
(\ref{measurepreserve1}) and (\ref{measurepreserve2}) give
\begin{eqnarray}
0 &=&\int_{\omega }B(\Theta )(h(\Theta )+s_{1})\,d\Theta ,
\label{m1} \\
0 &=&\left\{ \frac{N-1}{R}+\frac{A^{\prime }(R)}{A(R)}\right\} \cdot
\int_{\omega }B(\Theta )(h(\Theta )+s_{1})^{2}\,d\Theta +s_{2}\int_{\omega
}B(\Theta )\,d\Theta .
\label{m2}
\end{eqnarray}
Next assume that $\phi \in C^{2}(\Omega )$. Then, using (\ref{decomp}) and
the second condition in (\ref{diffPmu_1}) a short computation shows that
\begin{eqnarray*}
0 &\leq &\left\{ (N-2)(N-1)R^{N-3}A(R)+2(N-1)R^{N-2}(A^{\prime} (R))^{ N-1}A^{\prime
\prime }(R)\right\} \cdot  \\
&&\qquad \cdot \int_{\omega }B(\Theta )(h(\Theta )+s_{1})^{2}\,d\Theta  \\
&& + s_{2}\left\{ (N-1)R^{N-2}A(R)+R^{N-1}A^{\prime }(R)\right\} \int_{\omega
}B(\Theta )\,d\Theta  \\
&& + R^{N-3}A(R)\int_{\omega }B(\Theta )\left\vert \nabla _{\Theta }(h(\Theta
))+s_{1}\right\vert ^{2}\,d\Theta .
\end{eqnarray*}
Together with (\ref{m2}) this implies
\begin{eqnarray*}
0 &\leq &\left\{ -(N-1)R^{N-3}A(R)-R^{N-1}\frac{A^{\prime 2} (R)}{A(R)}
+R^{N-1}A^{\prime \prime }(R)\right\} \cdot  \\
&&\qquad \cdot \int_{\omega }B(\Theta )(h(\Theta )+s_{1})^{2}\,d\Theta  \\
&& + R^{N-3}A(R)\int_{\omega }B(\Theta )|\nabla _{\Theta }(h(\Theta
)+s_{1})|^{2}\,d\Theta .
\end{eqnarray*}
This implies
(\ref{stability}), in view of (\ref{m1}), and the definition of
$\lambda (B,\omega )$.
$\hfill \Box $
%\\[0.2cm]
\begin{remark}\label{Remark 2.3.} \rm  The value of $\lambda(B,\omega )$ is explicitly
known in some special cases. For instance (see, e.g. \cite{SW}
), if $B\equiv 1$, and
$\omega = \mathbb{S}^{N-1}$, we have
\begin{equation}
\label{eigensphere}
 \lambda (1,\mathbb{S}^{N-1})=N-1,
\end{equation}
the eigenvalue has multiplicity $N$, with corresponding
eigenfunctions $ u_{i}(x)=x_{i}$, ($i=1,\ldots ,N$), so that
(\ref{stability}) reads as
\begin{equation}
A^{\prime 2}\leq A^{\prime \prime }(r)A(r),  
\label{stabilityeasy}
\end{equation}
or equivalently, $A$ is log-convex, that is,
\begin{equation*}
A(r)=e^{g(r)},
\end{equation*}
with a convex function $g$. It has been conjectured in \cite{RCBM},
Conjecture 3.12, that for weights $\phi =A(r)$, with  log-convex
$A$, balls $B_{R}$, ($R>0$), solve the isoperimetric problem in
$\mathbb{R}^{N}$. 
\\
After finishing this paper, S. Howe kindly informed us about his new
preprint \cite{Howe} where he gives some partial answers to this
conjecture. He also determines the isoperimetric sets for some
radial weights. Further, some numerical evidence for the validity of the log-convex conjecture is provided in \cite{LMPW}.
\\
It is interesting to note that Theorem 1.1, whose proof will
be the object of the next section, and Theorem 2.1 imply the
following result.
\end{remark}
%\\[0.2cm]

\begin{proposition}\label{Proposition 2.4.} {\sl
Let $k\geq 0$, and
\begin{equation}
B=B_{k}(\Theta )=\left( \frac{x_{N}}{|x|}\right) ^{k},\quad (x\in
\mathbb{S}_{+}^{N-1}).
 \label{Bhalfsphere}
\end{equation}
Then
\begin{equation}
\lambda (B_{k},\mathbb{S}_{+}^{N-1}) = N-1+k,
\label{eigenk}
\end{equation}
with corresponding eigenfunctions
\begin{equation}
\label{eigenfunctionforxNk} u_i = x_i , \quad (i=1, \ldots , N-1).
\end{equation}
}
\end{proposition}

\noindent {\sl Proof: } Let $u_i $ be given by (\ref{eigenfunctionforxNk}).
Theorem 1.1 and Theorem 2.1 imply that
(\ref{stability}) holds, with $\omega = \mathbb{S}_+ ^{N-1} $,
$A(r) =r^k e^{cr^2 } $, ($c\geq 0$), and $B(\Theta ) = B_k (\Theta
)$. Hence $ \lambda ( B_k , \mathbb{S}_+ ^{N-1} ) \geq N-1 +k
-2cr^2 $ for all $r>0$, which implies that $\lambda ( B_k ,
\mathbb{S}_+ ^{N-1} ) \geq N-1 +k $. The assertion follows from the
identities
\begin{eqnarray*}
 & & \int_{\mathbb{S}_+ ^{N-1} } |\nabla _{\Theta } u_i |^2 B_k \, d\Theta
= (N-1+k) \int_{\mathbb{S}_+ ^{N-1} } (u_i )^2 B_k \, d\Theta  , 
\quad \mbox{and }\\
 & & \int_{\mathbb{S}_+ ^{N-1} } u_i  B_k \, d\Theta =0 ,
 \quad (i=1, \ldots , N-1 ) .\qquad \qquad \hfill \Box
\end{eqnarray*}

The next result gives the sharp constant in a weighted Hardy
inequality with respect to the measure $ x_N ^k |x| ^m \, dx$ in the
quarter space $\{ x_1 >0 , x_N >0 \} $
 (for related results in half spaces, see
e.g., \cite{AFV}, \cite{BFL}, \cite{Ma}, \cite{Nazarov} and \cite{Tid}).
\\
First we introduce some notation.  Let $D$ be an open set in
$\mathbb{R} ^N _+ $, and $\nu $ a measure given by $d\nu = \phi (x)dx
$, where  $\phi \in L_{\mbox{loc}} ^{\infty } (\mathbb{R} ^N _+ ) $, and
$\phi (x) > 0$.
The weighted H\"older space $L^2 (D , d\nu ) $ is the set of all measurable functions
$u: D \to \mathbb{R} $ such that
$ \int _D  u^2 \, d\nu   <+ \infty $,
and  
the weighted Sobolev space $W^{1,2} (D, d\nu ) $ is the set of functions 
$u\in L^2 (D, d\nu )$ that possess weak partial derivatives 
$u_{x_i } \in L^2 (D, d\nu )$, ($i=1, \ldots ,N$). 
Norms in these spaces are given respectively by
\begin{equation*}
\Vert u\Vert _{L^2 (D, d\nu )}  :=  \left( 
 \int _D  u^2 \, d\nu \right) ^{1/2} , 
 \end{equation*}
 and
\begin{equation*}
\Vert u\Vert _{W^{1,2} (D, d\nu )}  :=  \left( 
\int _D  (|u^2 +|\nabla  u|^2 )\, d\nu \right) ^{1/2} .
 \end{equation*}
\begin{definition}\label{Definition 2.5. } {\sl   
Let $X$ be the set of all functions 
$u \in C^1 (\overline{D} )$ that vanish in a neighborhood of 
$\partial D \setminus \{ x_N =0 \} $.  Then let 
$V^2 (D  ,d\nu )$ be the 
closure of $X$ in the norm of $W^{1,2} (D, d\nu )$.}
\end{definition}

%\\[0.2cm]
Next, let
\begin{equation}
\label{Q} Q:= \{ x\in \mathbb{R}^N :\, x_1
>0, \, x_N
>0 \} ,
\end{equation}
and specify
\begin{equation}
\label{wspecial} d \nu := x_N ^k |x|^m \, dx,
\end{equation}
where $k\geq 0 $ and $m\in \mathbb{N} $.
%\\[0.2cm]

\begin{theorem}\label{Theorem 2.6.} {\sl
 With $Q$ and $\nu $ given by (\ref{Q}) and
(\ref{wspecial}) respectively, we have
\begin{equation}
\label{hardyweight} \int _{Q } |\nabla u|^2  \, d\nu \geq C(k,m) \int
_{Q } \frac{u^2}{|x|^2} \, d\nu ,
\end{equation}
for all $u\in V^2 (Q,  d\nu )$,
 where
\begin{equation}
\label{hardyconstant} C(k,m) = \left( \frac{N+m+k -2}{2} \right) ^2
+ N+k-1 = \left( \frac{N+m+k}{2} \right) ^2 -m .
\end{equation}
The constant $C(k,m)$ in (\ref{hardyweight}) is sharp, and is not
attained for any nontrivial function $u$.
}
\end{theorem}

\noindent {\sl Proof : } We proceed as in [34, proof of Proposition
4.1].
\\
Extend $u$ to an odd function onto $\mathbb{R} ^N _+ $ by setting $
u(-x_1 , x_2 , \ldots ,x_N ) :=- u(x) $, ($ x\in Q$). Writing $u=
u(r, \Theta )$, and
 $B_k (\Theta ) = w(x)= x_N ^k |x|^{-k } $, we have for a.e. $r >0$,
 $$
 \int_{\mathbb{S}^{N-1} _+ } u (r, \Theta ) B_k (\Theta ) \,
 d\Theta =0,
 $$
 and thus
by \label{Proposition 2.4.},
\begin{equation}
\label{slice1} \int_{\mathbb{S}^{N-1} _+ } |\nabla _{\Theta } u (r,
\Theta )| ^2 B_k (\Theta ) \, d\Theta \geq (N+k -1 ) \int_{
\mathbb{S}^{N-1} _+ } [u  (r, \Theta )]^2 B_k (\Theta ) \, d\Theta .
\end{equation}
Further, the one-dimensional Hardy inequality (see \cite{Bliss})
tells us that for a.e. $\Theta \in \mathbb{S} ^{N-1} _+ $,
\begin{equation}
\label{slice2} \int_0 ^{+\infty } r^{N+m +k-1} [u_r (r, \Theta )] ^2
\, dr \geq \left( \frac{N+m+k -2}{2} \right) ^2 \int_0 ^{+\infty }
r^{N+m +k-3} [u (r, \Theta )]^2 \, dr .
\end{equation}
Integrating (\ref{slice1}) and (\ref{slice2}) gives
\begin{eqnarray*}
 \int _{\mathbb{R} ^N _+ } |\nabla u |^2  \, d\nu  
 & = &
\int _0 ^{+\infty } \int_{\mathbb{S}^{N-1} _+ } \left( [u_r ]^2 +
r^{-2} |\nabla _{\Theta } u | ^2 \right) r^{N-1 + m+k} B_k \, d\Theta \, dr
\\
 & \geq  &  \left[ \left( \frac{N+m+k -2}{2} \right) ^2
+ N+k-1 \right]   \int_0 ^{\infty } \int_{\mathbb{S}^{N-1} _+ } u^2
r^{N+ m+k-3} B_k \, d\Theta \, dr
\\
 & = & C(k,m) \int_{\mathbb{R}^N _+ } \frac{u^2}{|x|^2 }  \, d\nu  .
\end{eqnarray*}
The constant $C(k,m)$ is not attained since the constant is not
attained in the one-dimensional Hardy inequality. Moreover, the
exactness of $C(k,m)$ follows in a standard manner by considering
functions of the form $u= u_n = x_1 |x|^{(-N-m-k)/2 } \psi _n (|x|)
$, $(n\in \mathbb{N} )$, where
 $\psi _n \in C_0 ^{\infty } ((0, +\infty ))$ ,
$0\leq \psi _n \leq 1 $, $|\psi ' _n | \leq 4/n $,
 $\psi _n (t)= 0 $ for $t\in (0,(1/n)]\cup [2 n , +\infty ) $,
and $\psi _n (t) =1 $ for $t\in [(2/n),n]$, and then passing to
the limit $n\to \infty $.
The details are left to the reader.
$\hfill \Box $
\\[0.2cm]
Theorem 2.1 has some further consequences when the cone $\Omega $ contains the
wedge
\begin{equation*}
W_{+}:=
\{
x=(x_1 ,\ldots ,x_N ):\,x_i >0,\,i=1,\ldots ,N\},
\end{equation*}
and if
\begin{equation}
\phi (x)=\prod\limits_{i=1}^{N}\phi _{i}(x_{i}),
\label{prodcoord}
\end{equation}
for some smooth functions $\phi _i $, ($i=1,\ldots ,N$).
\\
In the following, let
\begin{equation*}
\omega _+ :=W_+ \cap \mathbb{S}^{N-1}.
\end{equation*}
We first show

\begin{lemma}\label{Lemma 2.7.} {\sl 
Assume that
$\phi \in C^2 (W_+)$ satisfies 
(\ref{decomp}) and (\ref{prodcoord}), 
where 
$ A,\phi_i
\in C^2 ((0,+\infty ))\cap C([0,+\infty ))$,
$B\in C^2(\omega _{+})\cap C(\overline{\omega _{+}})$, 
$\phi_i (x_i )>0$ for $ x_i >0$, 
($i=1,\ldots ,N$), 
$A(r)>0$ for
$r>0$, and $B(\Theta )>0$ for 
$ \Theta \in \omega _+$. Then
\begin{equation}
\phi(x)=a
\prod_{i=1}^N
x_i ^{k_i} e^{c|x|^{2}},
\quad x\in W_{+},
\label{decomp1}
\end{equation}
where $a>0$, $k_i \geq 0$, ($i=1,\ldots ,N$), and $c\in \mathbb{R}$.
}
\end{lemma}

%\\[0.2cm]
\noindent {\sl Proof : } Differentiating the equation $\log [ A(r)B(\Theta
)]=\log [\prod_{i=1}^{N}\phi _{i}(x_{i})]$ with respect to $r$ gives
\begin{equation*}
\frac{rA^{\prime }(r)}{A(r)}=\sum_{i=1}^{N}\frac{x_{i}\phi _{i}^{\prime
}(x_{i})}{\phi _{i}(x_{i})}. 
\end{equation*}
Differentiating this with respect to $x_i$ yields
\begin{equation*}
\frac{A^{\prime }(r)}{rA(r)}+ \frac{A^{\prime
\prime}(r)}{A(r)}-\frac{ (A^{\prime }  (r) )^2  }{(A (r))^2 } = \frac{\phi
_{i}^{\prime }(x_{i})}{x_{i} \phi _{i}(x_i )}+\frac{\phi _{i}^{\prime
\prime }(x_{i})}{\phi _{i}(x_{i})}- \frac{ (\phi _{i}^{\prime
}(x_{i}))^{2}}{(\phi _{i}(x_{i}))^{2}}=4c,\quad (i=1,\ldots ,N),
\end{equation*}
for some number $c\in \mathbb{R}$. 
In other words,
$$
\frac{d}{dx_i } \left\{ \frac{x_i \phi _i ^{\prime }(x_i)}{\phi _i (x_i )} \right\} = 
4cx_i  , \quad ( i=1, \ldots ,N) .
$$
Integrating this and dividing by $x_i $ give
$$
\frac{ \phi _i ^{\prime } (x_i )}{\phi _i (x_i ) } = 2cx_i + \frac{k_i }{x_i }, \quad (i=1, \ldots , N),
$$
for some numbers $k_i \in \mathbb{R} $, ($i= 1, \ldots ,N$).
Then another integration leads to   
\begin{equation*}
\log[ \phi _i (x_i)] =b_i + k_i \log x_i + c (x_i )^2 , \quad (b_i \in \mathbb{R}), 
\end{equation*}
that is,
$$
\phi _i (x_i ) = a_i x_i ^{k_i} e^{c(x_i )^2 },
$$
where $a_i = e^{b_i }$, ($i=1, \ldots , N$). Since $\phi _i \in C([0, +\infty ))$, and $\phi _i (x_i )>0$ for $x_i >0$, we have  $a_i >0$, and  $k_{i}\geq 0$, ($i=1, \ldots ,N$). Now  
(\ref{decomp1}) follows, (with $a= \prod_{i=1} ^N a_i $).
$\hfill \Box $
\\[0.2cm]
As pointed out in the Introduction,
we can specify the expression of the density $\phi$ of the measure,
when the cone $\Omega$ is $\mathbb{R}_{+}^{N}$ and  $\phi$  is factorized.
%\\[0.2cm]

\begin{theorem}\label{Theorem 2.8.} {\sl Assume $\Omega =\mathbb{R}_{+}^{N}$ and consider
Problem (\ref{Imu}), where $\phi \in C^{1}(\mathbb{R}_{+}^{N})\cap
C(\overline{\mathbb{R} _{+}^{N}})$, and satisfies (\ref{prodcoord}),
for some functions $\phi _{i}\in C^{2}(\mathbb{R})$, $\phi
_{i}(t)>0$ for $t\in \mathbb{R}$, ($i=1,\ldots ,N-1$), and $ \phi
_{N}\in C^{2}((0,+\infty ))\cap C([0,\infty ))$, $\phi _{N}(t)>0$
for $
t>0$. Suppose that $I_{\nu }(m)=P_{\nu }(B_{R}\cap \mathbb{R}_{+}^{N},
\mathbb{R}_{+}^{N})$ for $m=\nu (B_{R}\cap \mathbb{R}_{+}^{N})$. Then
\begin{equation}
\phi (x)=ax_{N}^{k}e^{c|x|^{2}},
\label{exp1}
\end{equation}
for some numbers $a>0$, $k\geq 0$ and $c\geq 0$.
}
\end{theorem}

\noindent {\sl Proof : } By Theorem 2.1 we have $\phi =A(r)B(\Theta )$
with smooth positive functions $A$ and $B$, and
\begin{equation}
\lambda (B,\mathbb{S}_{+}^{N-1})\geq N-1+r^{2}\left[ \frac{(A^{\prime
})^{2}}{A(r)^{2}}-\frac{A^{\prime \prime }(r)}{A(r)}\right] \quad \forall
r>0.  
\label{stability1}
\end{equation}
Then, Lemma \ref{Lemma 2.7.} shows that $\phi $ satisfies (\ref{decomp1}). 
Since $\varphi (x)>0$ whenever $x_N >0$, and $x_i =0 $, 
for some $i\in \{ 1, \ldots , N-1 \} $, it follows 
that we must have $k_i =0 $, ($i=1, \ldots , N-1 $). This proves  
(\ref{exp1}), for some numbers $a>0$,
$k\geq 0$ and $c \in \mathbb{R}$. Hence, $B(\Theta
)=[x_{N}|x|^{-1}]^{k}$, and $A(r)=ar^{k}e^{cr^{2}}$. Therefore (\ref
{stability1}) and (\ref{eigenk}) imply that
\begin{equation*}
N-1 +k\geq N-1 + k- 2cr ^2   \quad \forall r>0.
\end{equation*}
Hence we must have $c\geq 0$.
$\hfill \Box $
\\[0.2cm]
We end this section by analyzing the case where the cone
$\Omega$ is $\mathbb{R}^{N}\setminus \{0 \}$.

\begin{theorem}\label{Theorem 2.9.} {\sl 
Assume $\Omega =\mathbb{R}^{N}\setminus \{0 \}$ and
consider Problem (\ref{Imu}), with $\phi \in
C^{2}(\mathbb{R}^{N}\setminus \{0\})\cap C(\mathbb{R}^{N})$, $\phi
(x)>0$ for $x\not=0$, and satisfies (\ref {prodcoord}), where $\phi
_{i}\in C^{2}(\mathbb{R}\setminus \{0\})\cap C( \mathbb{R})$, and
$\phi _{i}(t)>0$ for $t\not=0$, ($i=1,\ldots ,N$). Suppose that
$I_{\nu }(m)=P_{\nu }(B_{R})$ for $m=\nu (B_{R})$. Then
\begin{equation}
\phi (x)=ae^{c|x|^{2}},
\label{exp2}
\end{equation}
for some numbers $a>0$, and $c\geq 0$.
}
\end{theorem}

\noindent {\sl Proof : } By Theorem 2.1 we have $\phi
=A(r)B(\Theta )$ with smooth positive functions $A$ and $B$, and
\begin{equation}
\lambda (B,\mathbb{S}^{N-1})\geq N-1+r^{2}\left[ \frac{A^{\prime }(r)^{2}
}{A(r)^{2}}-\frac{A^{\prime \prime }(r)}{A(r)}\right] \quad \forall r>0.
\label{stability2}
\end{equation}
Then, Lemma \ref{Lemma 2.7.} shows that $\phi $ satisfies (\ref{decomp1}).
Since $\varphi (x)>0$ whenever $x\not=0$ and $x_i =0 $, for some 
$i\in \{ 1, \ldots , N \} $, it follows that $k_i =0 $, ($i=1, \ldots , N $). 
This proves  
(\ref{exp2}), for some numbers  $
a>0$, and $c \in   \mathbb{R}$, that is, $B(\Theta )\equiv 1$ and
$A(r)=ae^{cr^{2}}$. Hence, (\ref{stability2}) and
(\ref{eigensphere}) imply that $A$ is log-convex, that is, we must
have $c\geq 0$.
$\hfill \Box $
\section{A Dido's problem}

In this section we provide the proof of Theorem 1.1. As pointed out
in the Introduction, we have to find the set having minimum $\mu$- perimeter
among all the subsets of $\mathbb{R}_{+}^{N}$ having prescribed $\mu$ -
measure, where $\mu$ is the measure defined in \eqref{intro_dmu}. In
order to face such a problem we first show a simple inequality for
measures defined on the real line related to $d\mu $. Then the
isoperimetric problem is addressed in the plane: the one-dimensional
results allow to restrict the search of optimal sets to the ones
which are starlike with respect to the origin. Finally Theorem
1.1 is achieved in its full generality.
\subsection{Dido's problem on the real line}

Let $\mathbb{R}_{+}=(0,+\infty ).$ The following isoperimetric inequality
holds.

\begin{proposition}\label{Proposition 3.1.}
{\sl Let $\phi :\mathbb{R}_{+}\rightarrow \mathbb{R}_{+}$ be a
nondecreasing continuous function, $d\nu =\phi (x)dx$ and $M$ be a
measurable subset of $\mathbb{R}_{+}$ with $\nu (M) <+ \infty $. Then
\begin{equation}
P_{\nu }(M)\geq P_{\nu }(S(M)),  
\label{1.1}
\end{equation}
where $S(M)$ denotes the interval $(0,d)$, with $d\geq 0$ chosen such that $
\nu (M)=\nu (S(M))$. }
\end{proposition}

\noindent {\sl Proof : } First assume that $M$ is of the form
\begin{equation}
\label{smoothsets}
M=\cup _{j=1}^{k}\left( a_{j},b_{j}\right) ,
\end{equation}
with
\begin{equation*}
0\leq a_{j}<a_{j+1},\text{ \ }a_{j}<b_{j},\text{ \ }b_{j}<b_{j+1}<+\infty ,
\end{equation*}
for all $j\in \left\{ 1,...,k-1\right\} $.
By the properties of the weight function $\phi $ we have that $b_{k}\geq
d $ and hence
\begin{equation}
\label{perimeterforsmoothsets}
P_{\nu }(M)=\sum\limits_{j=1}^{k}\left[ \phi (a_{j})+\phi (b_{j})
\right] \geq \phi (0)+\phi (d)=P_{\nu }(S(M)).
\end{equation}
Next let $M$ be measurable and $\nu (M)<+ \infty $.
By the basic properties of the perimeter, there exists a sequence of sets $\{ M_n \} $ of the form (\ref{smoothsets}) such that $\lim_{n\to +\infty } \nu (M\Delta M_n ) =0$ and 
$\lim_{n\to +\infty }  P_{\nu } (M_n )= P_{\nu } (M)$. The first limit implies that also
$\lim_{n\to +\infty } P_{\nu } ( S(M_n )  ) = P_{\nu } (S(M))$, so that the assertion follows from inequality (\ref{perimeterforsmoothsets}).
$\hfill \Box $
\subsection{Dido's problem in two dimensions}
In our study of the measure $d\mu ,$ an important role will be played by the
following isoperimetric theorem (see \cite{BMP} and \cite{RCBM}) relative to
the measure
\begin{equation*}
d\tau =\exp ( c\left\vert x\right\vert ^2 ) dx,\text{ \ \ }x\in
\mathbb{R}^{m},\text{ with }m\geq 1\text{ and }c\geq 0.
\end{equation*}

\begin{theorem}\label{Theorem 3.2. } 
{\sl If $M$ is any measurable subset of $\mathbb{R}^{N}$
and $ M^{\star }$ is the ball of $\mathbb{R}^{N}$ centered
at the origin having the same $\tau -$measure of $M$, then
\begin{equation}
P_{\tau }(M)\geq P_{\tau }(M^{\star  }).
\end{equation}
}
\end{theorem}

We write $(x,y)$ for points in $\mathbb{R}^{2}$, and we consider in $\mathbb{R}
_{+}^{2}$ the measure
\begin{equation*}
d\mu =y^{k}\exp \left( c(x^{2}+y^{2})\right) \, dx\, dy,
\end{equation*}
where $c\geq  0$ and $k\geq 0$. If $M$ is a measurable subset of $\mathbb{R}
_{+}^{2}$, given any number $m>0$, the isoperimetric problem on $\mathbb{R}
_{+}^{2}$ reads as:
\begin{equation}
I_{\mu }(m):=\inf \{P_{\mu }(M),\,\text{with }M:\mu (M)=m\}.  
\label{isop1}
\end{equation}
The following result holds true.

\begin{theorem}\label{Theorem 3.3. } {\sl Let $m>0$. Then $I_{\mu }(m)$ is attained for the
half-disk $B_{r}\cap \mathbb{R}_{+}^{2}$, centered at zero, having
$\mu $-measure $m$. Equivalently there exists $r>0$ such that
\begin{equation}
I_{\mu }(m)=P_{\mu }(B_{r}\cap \mathbb{R}_{+}^{2})=\exp \left(
cr^{2}\right) r^{k+1}\int_{0}^{\pi }\sin ^{k}\theta \, d\theta =
B\left( \frac{k+1}{2} , \frac{1}{2} \right) \exp \left(
cr^{2}\right) r^{k+1},  \label{halfcircle}
\end{equation}
where $B$ denotes the Beta function.
}
\end{theorem}

\noindent {\sl Proof: } If $k=0$, and $c=0$ (unweighted case), the result is well-known. Further, if $c>0$ and $k=0$, that is, $d\mu = e ^{c(x^2 + y ^2 )} \, dx\, dy $, the result follows from Theorem \ref{Theorem 3.2. } via reflexion about the $x$-axis.
Finally, the result has been shown in the case $c=0$ and $k>0$ by Maderna and Salsa,
\cite{MS}, (see also \cite{EMMP}).
\\
Therefore we may restrict ourselves to the case that both $c$ and $k$ are positive.
 
Our proof requires some technical effort which is mainly due to the degeneracy of the measure on the $x$-axis.
The strategy is as follows: 
First we use symmetrization arguments in order to reduce the isoperimetric problem to sets which are starlike w.r.t. the origin ({\sl Step 1}).
Then we obtain some a-priori-estimates for a minimizing sequence ({\sl Step 2}). This   allows us to show that a (starlike) minimizer exists ({\sl Step 3}), which is also bounded ({\sl Step 4}) and smooth ({\sl Step 5}). In {\sl Step 6 } we evaluate the second variation of the Perimeter functional, and we show that the minimizer is a half-disk centered at the origin.  

Throughout our proof, $C$ will denote a generic constant which may vary from line to line.   
\\[0.1cm]
{\sl Step 1: Symmetrization } 
\\[0.1cm]
Our aim is to simplify the isoperimetric problem using Steiner symmetrization in two directions. This method has already been employed in the case $c=0$ (see \cite{MS}, and
\cite{EMMP}).
\\     
Let $\left\{ D_{n}\right\} \subset \mathbb{R}_{+}^{2}$ be a
minimizing sequence for problem (\ref{isop1}), i.e.
\begin{equation*}
\mu (D_{n})=m \ \  \forall n\in \mathbb{N} \ \mbox{  and } \ \lim_{n\to
+\infty }P_{\mu }(D_{n})=I_{\mu }(m),
\end{equation*}
where, without loss of generality, we may assume that the sets $D_{n}$ are
smooth.

Let $D$ be a smooth set of $\mathbb{R}_+ ^2 $. We denote by $S_x (D)$
and $ S_y (D)$ the Steiner symmetrization in $x$-direction, with
respect to the measure $d\mu _{x}= e^{cx^2 } dx$, and the Steiner symmetrization in $y$-direction, with respect to
the measure $d\mu _y =e^{cy^2 }y^k \, dy$, of $D$, respectively.

More precisely, $S_x (D)$ is the subset of $\mathbb{R}_+ ^2$
whose cross sections parallel to the $x$-axis are open intervals
centered at the $y$-axis, and such that their $\mu _{x}$-lengths are equal to
those of the corresponding cross sections of $D$.

The set $S_y (D)$ is defined in a similar way: its cross
sections parallel to the $y$-axis are open intervals with an endpoint lying
on the $x$-axis, and such that their $\mu _{y}$-lengths are equal to those of the
corresponding cross sections of $D$.

Now consider the sequence of sets $M_{n}=S_y (S_x(D_n))$. By
Proposition \ref{Proposition 3.1.} and Theorem \ref{Theorem 3.2. } we have that
$P_{\mu }(S_y (S_x (D_n )))\leq P_{\mu }(D_n )$ and by Cavalieri's
principle $\mu (S_y (S_x (D_n )))=\mu (D_n )$. 
Therefore $\{ M_n \} $
is still a minimizing sequence for (\ref{isop1}). On one hand, the
sets $M_n $ can lose regularity under symmetrization: the
symmetrized sets are not more then locally Lipschitz continuous, in
general. On the other hand, they acquire some nice geometrical
property: they are all starlike with respect to the origin.
Thus, introducing polar coordinates $(r,\theta )$ by $x=r\cos \theta
$ and $y=r\sin \theta $, we have
\begin{equation}
M_{n}=\{(r,\theta ):\,0< r<\rho _{n}(\theta ),\,\theta \in (
0,\pi ) \},\quad \forall n\in \mathbb{N}, 
\label{starlike}
\end{equation}
for some functions $\rho _{n}(\theta ):( 0,\pi ) \rightarrow \left(
0,+\infty \right) .$ Note that, defining $\rho _n (0):= \lim
_{\theta \to 0 ^+ } \rho _n (\theta )=: \rho _n (\pi ) $, and 
$ \rho _n (\pi /2 ) :=
\lim_{\theta \to \pi /2 ^-} \rho _n (\theta )$, we have also have
$\rho _n \in C([0,\pi ] )$. Then
\begin{itemize}
\item[(i)] the functions $\rho _{n}(\theta )$ are locally Lipschitz in $
\left( 0,\pi /2 \right)$;
\item[(ii)] $\rho _{n}(\theta )=\rho _{n}(\pi -\theta ),\ \forall n\in \mathbb{
N},\ \forall \theta \in (0,\pi ) $;
\item[(iii)] the functions $x_{n}(\theta ):= \rho _{n}(\theta ) \cos \theta $ and $
y_{n}(\theta ):=\rho _{n}(\theta ) \sin \theta $ are
nonincreasing and nondecreasing, respectively, on $(0,\pi /2)$.
\end{itemize}
Hence we may assume that the minimizing sequence is of the form
(\ref{starlike}), with conditions (i)--(iii) in force. Under these
conditions, the set $M_{n},$ its $\mu $-measure and $\mu $-perimeter
are uniquely determined by
the function $\rho _{n}(\theta )$. More precisely, setting 
\begin{eqnarray*}
z & := & \sin ^{k}\theta, \quad 
\theta \in [ 0,\pi ],
\\
F(r) & := & \int_{0}^{r} e^{cr^2 }  t^{k+1}\,
dt, \quad
\mbox{and}
\\
G(r,p) & := & e^{cr^2 }  r^k \sqrt{r^{2}+p^{2}},\quad r>0,\
p\in \mathbb{R},
\end{eqnarray*}
we find that
\begin{eqnarray*}
\mu (M_{n}) & = & \int_{0}^{\pi }F(\rho _{n})z\,
d\theta =:\mu (\rho _{n}),
\quad \mbox{and }\\
P_{\mu }(M_{n}) & = & \int_{0}^{\pi }G(\rho _{n},\rho _{n}^{\prime })z\,
d\theta
=:P_{\mu }(\rho _{n}).
\end{eqnarray*}
With this notation, the isoperimetric problem (\ref{isop1}) now
reads as
\begin{eqnarray}
\label{isop2}
 \ \ \mbox{Minimize} & & P_{\mu }(\rho )  \ \mbox{over} 
\\
\nonumber 
 & & K  :=  \{\rho :(0,\pi /2
)\cup (\pi /2 , \pi )\to (0,+\infty ): \, \rho \mbox{ 
satisfies (i)-(iii) and } \mu (\rho )=m \}. 
\end{eqnarray}
{\sl Step 2: Some estimates }
\\[0.1cm]
Next we will obtain some uniform estimates for the minimizing sequence
$\{\rho _{n}\} $ of problem (\ref{isop2}). 
 
Condition (iii) implies  
\begin{equation}
-\rho _{n}(\theta )\cot \theta
\leq \rho _n ^{\prime }(\theta )\leq \rho _n (\theta )
\tan \theta \ \ \ \mbox{ a.e. on }\ (0,\pi /2),\
n\in \mathbb{N}
.  
\label{cond2}
\end{equation}
Set
\begin{equation*}
y_n ^0 :=\sup_{\theta \in (0,\pi /2)} y_n (\theta )=
y_n (\pi /2) = \rho _n (\pi /2 ).
\end{equation*}
We claim that
\begin{equation}
\sup_{n\in \mathbb{N}}
y_n ^0 =:
y^0 <+\infty .
\label{cond3}
\end{equation}
Indeed, since $\{P_{\mu } (\rho _n ) \} $ is a bounded sequence, 
we obtain for every $n\in \mathbb{N}$,
\begin{eqnarray*}
C & \geq & P_{\mu }(\rho _n )
\\ 
 & = &
2\int_0 ^{\pi /2 }
e^{  c(
x_n ^2 (\theta )+y_n ^2 (\theta )) }
 y_n ^{k}
(\theta )
\sqrt{
(x_n ^{\prime }(\theta ))^2 +(y_n ^{\prime }(\theta ))^2 } \, 
d\theta \\
& \geq &
2\int_0 ^{\pi /2}e^{  cy_n ^2 (\theta )}
y_n ^k (\theta ) y_n ^{\prime }(\theta )d\theta = 2\int_0 ^{y_n ^0 }
e^{
 ct^2 } t^k dt ,
\end{eqnarray*}
and (\ref{cond3}) follows.

From (\ref{cond2}) and (\ref{cond3}) we further deduce that 
for every $\theta
 \in (0,\pi )$, 
\begin{equation}
\rho _n (\theta )=\frac{y_n (\theta )}{\sin \theta }\leq \frac{y_n (\pi
/2)}{\sin \theta }\leq \frac{y^0 }{\sin \theta } \quad \forall n\in
\mathbb{N}.
\label{rho_n<C/sin}
\end{equation}

Conditions 
(\ref{rho_n<C/sin}) and (\ref{cond2}) imply that for every 
$\delta \in (0, \pi /4)$ there is a number $d_{\delta } >0 $ such that 
\begin{equation}
\sup_{\theta \in (\delta ,\pi /2-\delta )}\left\{ \rho _{n}(\theta
),\left\vert \rho _{n}^{\prime }(\theta )\right\vert \right\} \leq d_{\delta
}.  
\label{A3}
\end{equation}

Next we claim:
\\
There exists a number $d_1 >0 $, such that 
\begin{equation}
\label{apriori2}
\rho _n (\theta ) \geq d_1  \quad \forall \theta \in (0, \pi ), \ \mbox{ and } \forall n\in \mathbb{N} .
\end{equation}
Assume (\ref{apriori2}) was not true. Then the fact that $ x_n (\theta ) $ and $ y_n (\theta )  $ are nonincreasing, respectively nondecreasing, $\forall n\in \mathbb{N} $, means that  
there is a subsequence, still labelled as $\{\rho _n \} $, 
such that $\lim _{n\to \infty } \rho _n (\pi /4 )  = 0 $. Set $\delta _n := \rho _n (\pi /4) /\sqrt{2} $, and note that $x_n (\pi /4 ) = y_n (\pi /4 ) = \delta _n $. 
In view of (\ref{rho_n<C/sin})  we have that 
$$
\lim _{n\to \infty } \mu (M_n \cap \{ |x|< \delta _n \} ) =0.
$$
Since $\mu (\rho _n )=m$, this implies that there is a number $d_2 >0 $, such that for all $n\in \mathbb{N}$, 
\begin{eqnarray}
\label{lim1}
d_2 & \leq &
\mu (M_n \cap \{ x>\delta _n \}) = - \int_0 ^{\pi /4} x_n ^{\prime } (\theta ) 
e^{cx_n ^2 (\theta )  } \int_0 ^{y_n (\theta )} e^{ct^2 } t^{k+1} \, dt \, d\theta  
\\
\nonumber
 & \leq & - (\delta _n ) ^2 \int_0 ^{\pi /4}  x_n ^{\prime } (\theta )
 e^{c(x_n ^2 (\theta ) + y_n ^2 (\theta ) )} y_n ^k (\theta  ) \, d\theta .
\end{eqnarray}  
On the other hand, the sequence $\{ P_{\mu } (\rho _n ) \} $ is bounded, so that 
\begin{eqnarray}
\label{lim2} 
C & \geq & 
P_{\mu } (\rho _n ) \\
\nonumber  & \geq & \int_0 ^{\pi /4}  
e^{c (x_n ^2 (\theta ) + y_n ^2 (\theta )  )} y _n ^k (\theta )
 \sqrt{(x_n ^{\prime } (\theta ))^2    + (y_n ^{\prime } (\theta ) ) ^2 } \, d\theta  
 \\
 \nonumber
  & \geq &  - \int_0 ^{\pi /4} x_n ^{\prime } (\theta ) 
  e ^{c(x_n ^2 (\theta ) + y_n ^2 (\theta ))} y_n ^k (\theta) 
 \, d\theta .
\end{eqnarray}
Hence we obtain $d_2 \leq \delta _n ^2 C$ for all $n\in \mathbb{N}$, which is a contradiction.

Next we claim that there is a number $d_3 >0$ such that holds for every 
$\theta \in (0, \pi /2 )$ and for all $n\in \mathbb{N}$,  
\begin{equation}
y_n ^k (\theta )\int_0 ^{x_n (\theta )} e ^{ct^2 } \,
dt
\leq d_3 .
\label{cond4}
\end{equation}
Consider the set
\begin{equation*}
\widetilde{M}_n (\theta) :=
\left\{ (x,y)\in M_n :
y\leq y_n (\theta )\right\} .
\end{equation*}
It is easy to verify that
\begin{eqnarray*}
\frac{1}{2}\left( P_{\mu }(M_n )-P_{\mu }
(\widetilde{M}_n (\theta))
\right) 
 & = & \int_{\theta }^{\pi /2} e ^{ c( x_n ^2 (\tau
)+y_n ^2 (\tau ))}  y_n (\tau
)^k \sqrt{(x_n ^{\prime }(\tau ))^2 +(y_n ^{\prime }(\tau ))^2 }\, d\tau 
\\
 & & 
-\int_0 ^{x_n (\theta )} e^{ c ( t^{2}+y_n ^2 (\theta
))}  y_n ^k (\theta )\, dt
\\
 & \geq & \int_{\theta }^{\pi /2} (-x_n ^{\prime } (\theta ) ) e^{  c( x_n ^2 (\tau
)}
\left( e^{y_n ^2 (\tau ) } y_{n}^{k}(\tau ) - e ^{y_n ^2 (\theta )} 
y_n ^k (\theta )
\right) \,  d\tau 
  \geq  0.
\end{eqnarray*}
Hence
\begin{equation}
\label{int55}
C\geq P_{\mu }(M_n )\geq P_{\mu }(\widetilde{M}_n (\theta ))\geq
2y_n ^k (\theta )e^{cy_n ^2 (\theta )}
\int_0 ^{x_n (\theta )} e^{ct^2 } \,  dt,
\end{equation}
and (\ref{cond4}) follows.

Below
we will frequently make use of the following limit which holds for all $\alpha >-1$,
\begin{equation}
\label{explimit}
\lim_{z\to +\infty } \frac{ \int_0 ^z e^{ct^2 } t^{\alpha } 
\, dt}{e^{cz^2 } z^{\alpha -1} } = \frac{1}{2c} .
\end{equation}
In view of (\ref{explimit}) with $\alpha =0$, and (\ref{int55}), and since $x_n (\theta )\geq d/\sqrt{2}$ for $\theta \in (0, \pi /4)$, 
we obtain
\begin{equation} 
\label{56}
C \geq y_n ^k (\theta ) e^{cy_n ^2 (\theta ) } \frac{e ^{cx_n ^2 (\theta )}}{ x_n (\theta )}, \quad \forall \theta \in (0, \pi /4) .
\end{equation}
Since   $y_n (\theta ) \geq (1/2 ) \theta \rho _n (\theta ) $ for $\theta \in (0, \pi /4 )$, and 
$ x_n (\theta ) \leq \rho _n (\theta )$, we further deduce from (\ref{56}),
\begin{equation}
\label{anotherestimate}
C\geq \rho _n ^{k-1} (\theta ) \theta ^k e^{c\rho _n ^2 (\theta )}, \quad  \forall \theta \in (0, \pi /4).
\end{equation}
Now recall (\ref{apriori2}), and $\lim _{z\to + \infty } e^{cz^ 2 /2 } z^{k-1} 
=+\infty $. Hence (\ref{anotherestimate}) shows 
that there is a number $d_4 >0 $ such that for all $n \in \mathbb{N}$,
\begin{equation}
\label{apriori3}
\rho_n (\theta ) \leq \sqrt{d_4  - \frac{2k}{c} \ln \theta } , \quad 
\forall \theta \in (0, \pi /4 ).
\end{equation}
Finally we show:
\begin{eqnarray}
\label{intestimate}
 & & \mbox{For every $\epsilon \in (0, m ) $ there is a 
 $\delta \in (0, \pi /2 )$, such that}
\\
\nonumber
 & &  \mu (M_n \cap \{ \delta <\theta <\pi -\delta  \} )  > m-  
 \epsilon, \quad \forall n\in \mathbb{N}.
\end{eqnarray}
Indeed, (\ref{anotherestimate}), (\ref{apriori3}) and (\ref{explimit}) 
with $\alpha =k+1$, show that 
\begin{eqnarray}
\label{sest}
\mu (M \cap \{  0<\theta < s \} )
 & = & 
 \int _0 ^s \sin ^k \theta \int _0 ^{\rho _n (\theta )}  
 e^{ct^2 } t^{k+1} \, dt \, d\theta 
\\
\nonumber
 & \leq & C \int_0 ^s  \theta ^k e^{c\rho _n ^2 (\theta ) } \rho _n ^k 
 (\theta ) \, d\theta   
\\
\nonumber
 & \leq & C \int_0 ^s   \sqrt{ d_4 - \frac{2k}{c} \ln \theta } \, d\theta 
 \rightarrow 0, \quad \mbox{ as } \ s\to 0.  
\end{eqnarray}
Now the claim (\ref{intestimate}) follows from the uniform estimate (\ref{sest}) and from the fact that
for every $s\in (0, \pi /2 )$, 
$$
m/2 = \mu (M_n \cap \{0< \theta <s \} )+ \mu (M_n \cap \{ s < \theta < \pi /2\} ) .
$$
{\sl Step 3: The minimum is achieved}
\\[0.1cm]
In this step we show that a minimizer of problem (\ref{isop2}) exists.

In view of the properties (i)--(iii), (\ref{cond2}),  and the estimates (\ref{rho_n<C/sin}), (\ref{apriori2}), (\ref{A3}) and (\ref{apriori3})
there exists a function $\rho : (0, \pi /2 ) \cup (\pi /2, \pi ) \to [0,+\infty )$ which is locally Lipschitz continuous, and a subsequence, still denoted by 
$\{ \rho _n \} $, such 
\begin{eqnarray}
\rho _ n & \rightarrow & \rho \text{ uniformly on compact subsets of
}(0,\pi
/2),  
\label{cond8} \\
\rho (\theta ) & = & \rho (\pi -\theta )\text{ } \ \forall \theta
\in (0,\pi /2),
\label{cond5} 
\\
-\rho (\theta )\cot \theta & \leq & \rho ^{\prime }(\theta )\leq
\rho (\theta )\tan \theta \quad \mbox{ a.e. on } \ (0,\pi /2),
\label{cond6}
\\
\label{cond10}
\rho (\theta ) & \leq & \frac{y^{0}}{\sin \theta }\quad 
\forall \theta \in (0, \pi /2),
\\
\rho (\theta ) & \geq & d_1 \quad \forall \theta \in (0, \pi /2 ), 
\label{cond11}
\\
 \sup_{\theta \in (\delta ,\pi /2-\delta )}\left\{ \rho (\theta
),\left\vert \rho ^{\prime }(\theta )\right\vert \right\} & \leq & d_{\delta
}, \quad \forall \delta \in (0, \pi /4 ),
\label{cond12}
\\
\rho (\theta ) & \leq & \sqrt{d_4  - \frac{2k}{c} \ln \theta } , 
\quad \forall \theta \in (0, \pi /4 ).
\label{cond13}   
\end{eqnarray}
Note, setting $x(\theta ):= \rho (\theta ) \cos \theta $, and $y(\theta ):= \rho (\theta ) \sin \theta $, condition (\ref{cond6}) implies that the functions $x(\theta )$ and $y(\theta )$ are nonincreasing, respectively nondecreasing on $(0, \pi /2)$.
Further, defining
$\rho (\pi /2 ):= \lim _{\theta \to \pi /2} \rho (\theta ) $, we see that
$\rho \in C((0,\pi ))$. 

Let 
$$
M:= \{ (r, \theta ): \ 0< r< \rho (\theta ),\ \theta \in (0, \pi )\} ,
$$ 
and $P_{\mu } (\rho ):= P_{\mu } (M)$.
We claim   
\begin{equation}
\label{fullmeasure}
\mu (M)=m.
\end{equation}
Indeed, the estimate (\ref{intestimate})  shows
\begin{eqnarray*}
 & & \mbox{For every $\epsilon \in (0, m)$ there is a 
 $\delta \in (0, \pi /2 )$, such that}
\\
\nonumber
 & &  \mu (M \cap \{ \delta <\theta <\pi -\delta   \} ) \geq m - \epsilon .
\end{eqnarray*}
Since we also have $\mu (M) \leq m$, (\ref{fullmeasure}) follows.
\\
Finally, the 
lower semicontinuity of the perimeter shows that
\begin{equation}
\label{semicontperimeter}
I_{\mu } (m) = \lim _{n\to \infty } P_{\mu } (M_n ) \geq P_{\mu } (M).
 \end{equation}
 But $\rho \in K$, therefore $I_{\mu } (m) = P_{\mu } (M)$, and $M$ is a minimizer. 
 
Note that our weight function $ \phi (x,y) := y^k e ^{c(x^2 + y^2 )} $ is positive and $\phi \in C^{\infty }( \mathbb{R} ^2 _+) $. 
Due to a regularity result of F. Morgan, \cite{Mo1}, Corollary 3.7 and Remark 3.10, this implies that $\partial M \cap \mathbb{R}^2 _+ $ is a one-dimensional
$C^1$-manifold which is locally analytic. In view of the symmetry of $M$ this implies that $\rho $ is differentiable at     
$\pi /2 $, with $\lim _{\theta \to \pi /2 } \rho ^{\prime } (\theta ) = \rho ^{\prime } (\pi /2)=0$. 
Using the properties (\ref{cond5})--(\ref{cond13}) this implies that $\rho \in C^{\infty } ((0, \pi ))$. 
\\ 
Then standard Calculus of Variations  (see \cite{MS})  shows that there is a number $
\gamma \in \mathbb{R}$ - a Lagrangian multiplier - such that  
\begin{equation}
-\frac{d}{d\theta }   \left( G_p z \right)   + G_r z=  \gamma
 F^{\prime } z \quad \mbox{on }\ (0, \pi ) .
\label{Euler}
\end{equation} 
Here and in the following, the functions $G,F$ and 
their derivatives are evaluated at $(\rho , \rho ^{\prime })$.
\\[0.2cm]
{\sl Step 4 : The minimizer is bounded }
\\[0.1cm]
We will argue by contradiction, that is, we assume 
that $\rho $ was unbounded. 
Then (\ref{cond6}) would imply that
\begin{equation}
\label{wronglimit} 
\lim_{t\to 0} \rho (t) =+\infty .
\end{equation}
First we claim that (\ref{wronglimit}) further means  
that there exists a sequence $t_n \to 0$ such that 
\begin{equation}
\label{rhoprimeest}
- \rho ^{\prime } (t_n ) \geq \rho ^3 (t_n ) .
\end{equation}
Indeed,
assume (\ref{rhoprimeest}) was not true. 
Then there exists a number $t_0 >0$ such that 
\begin{equation}
\label{rhoprimeest2}
-\rho ^{\prime } (t) < \rho ^3 (t) \quad \mbox{for $t\in (0, t_0 )$.}
\end{equation} 
By the estimate (\ref{cond13}) we can find a number 
$t_1 \in (0, t_0 )$ such that
$-2t_1 + (\rho  (t_1 ))^{-2} =: \delta _0 >0 $.
Integrating (\ref{rhoprimeest2}) gives
$$
\frac{1}{\rho ^2 (t) }  >  2t -2 t _1 + \frac{1}{(\rho  (t_1 ))^2 } = 
\delta _0 + 2t 
\quad \forall t\in (0, t_1 ),
$$
which implies that $\rho $ is bounded, a contradiction. Hence 
(\ref{rhoprimeest}) follows.
Note that (\ref{rhoprimeest}), together with our assumption 
(\ref{wronglimit}) implies that 
\begin{equation}
\label{rhorhoprime}
\lim_{n\to \infty } \rho (t_n )/ \rho ^{\prime } (t_n )=0.
\end{equation} 
Using the Euler equation 
(\ref{Euler}), a short calculation shows that
\begin{equation}
\label{Eulermodified}
\frac{d}{d\theta } \left(
G - \rho ^{\prime } G_p   - \gamma F \right)  = \rho ^{\prime } G_p 
\frac{z^{\prime}}{z} .
\end{equation}
Integrating (\ref{Eulermodified}) on the interval  $(t_n , \pi /2 )$
gives
\begin{eqnarray}
\label{integraleuler}
 & & \gamma \int_0 ^{\rho (t_n )} e^{cs^2 } s^{k+1 } ds - 
 e^{c(\rho (t_n ))^2} (\rho (t_n ))^{k+2} \left( (\rho (t_n ))^2 + 
 (\rho ^{\prime} (t_n ))^2 \right) ^{-1/2} 
 \\
\nonumber 
 & = & -c_1 +  \int_{t_n } ^{\pi /2 }  e^{c(\rho (t))^2 } 
 (\rho (t))^{k}  (\rho ^{\prime } (t)) ^2 
 \left( (\rho (t ))^2 + (\rho ^{\prime} (t))^2 
 \right) ^{-1/2} k \cot t \, dt ,
\end{eqnarray} 
where we have put
$$
c_1 = (G -\rho ' G_p -\gamma F ) \Big| _{\theta = \pi /2 }.
$$       
In view of (\ref{wronglimit}), 
(\ref{rhoprimeest}), (\ref{rhorhoprime}) 
and (\ref{explimit}) with $\alpha= k+1$ we find that
\begin{equation}
\label{ratio1}
\lim_{n\to \infty } \frac{\gamma \int_0 ^{\rho (t_n )} 
e^{cs^2 } s^{k+1 } ds}{e^{c(\rho (t_n ))^2} (\rho (t_n ))^{k+2} \left( 
(\rho (t_n ))^2 + (\rho ^{\prime} (t_n ))^2 \right) ^{-1/2} } =+\infty .
\end{equation}
Hence the left-hand side of equation (\ref{integraleuler}) 
tends to $+\infty $ as $n \to + \infty $. 
Using de l'Hospital's rule, (\ref{wronglimit}), (\ref{rhoprimeest}) and (\ref{rhorhoprime}), we obtain from (\ref{ratio1}),
\begin{eqnarray*}
1 & = & \lim_{n\to \infty }
\frac{
\gamma \int_0 ^{\rho (t_n )} e^{cs^2 } s^{k+1 } ds - e^{c(\rho (t_n ))^2} (\rho (t_n ))^{k+2} \left( (\rho (t_n ))^2 + (\rho ^{\prime} (t_n ))^2 \right) ^{-1/2}
 }{
 -c_1 +  \int_{t_n } ^{\pi /2 }  e^{c(\rho (t))^2 } (\rho (t))^k  
 (\rho ^{\prime } (t)) ^2  
 \left( (\rho (t ))^2 + (\rho ^{\prime} (t))^2 \right) ^{-1/2}  
 k \cot t \, dt}
\\
 & = & 
\lim_{n\to \infty }
\frac{
\gamma \int_0 ^{\rho (t_n )} e^{cs^2 } s^{k+1 } ds 
 }{ 
 \int_{t_n } ^{\pi /2 }  e^{c(\rho (t))^2 } (\rho (t))^k    
 (\rho ^{\prime } (t)) ^2 
 \left( (\rho (t ))^2 + (\rho ^{\prime} (t))^2 \right) ^{-1/2} 
 k \cot t \, dt} 
\\
 & = & \lim_{n\to \infty }
\frac{\gamma \rho ^{\prime } (t_n )  e^{c(\rho (t_n ))^2 } 
(\rho (t_n ))^{k+1 }  
 }{ 
-  e^{c(\rho (t_n ))^2 } (\rho (t_n ))^k   
( \rho ^{\prime} (t_n ) )^2  \left( (\rho (t_n ))^2 + (\rho ^{\prime} (t_n ))^2 
\right) ^{-1/2}
k \cot t_n }
\\
 & = & \lim_{n\to \infty }
\frac{\gamma \rho ^{\prime } (t_n )  
e^{c(\rho (t_n ))^2 } (\rho (t_n ))^{k+1 }  
 }{ 
  e^{c(\rho (t_n ))^2 } (\rho (t_n ))^k   
 \rho ^{\prime} (t_n )  
k \cot t_n }
\\
 & = &
  \lim_{n\to \infty } 
  \frac{
  \gamma \rho (t_n )
  }{ 
  k\cot t_n }
 =  \lim_{n\to \infty } 
 \frac{\gamma }{k} t_n \rho (t_n ). 
\end{eqnarray*} 
But the last limit is zero in view of (\ref{cond13}), and we have obtained a contradiction. In other words,  
$ \rho $ is bounded on $(0, \pi )$.

Putting $\rho (0) := \lim _{t\to 0} \rho (t) =: \rho (\pi )$, we then have 
\begin{equation}
\label{rhocont} 
\rho \in C([0, \pi ]).
\end{equation} 
{\sl Step 5: $\rho ^{\prime \prime } $ is bounded }
\\[0.1cm]
We will first need some integrability properties of the functions
\begin{eqnarray*}
G_r & = & e^{c\rho ^2 } \rho ^k \left( [2c\rho + (k/\rho )]  \{ \rho ^2 + (\rho ^{\prime } ) ^2 \} ^{1/2}  + \rho  \{ \rho ^2 + (\rho ^{\prime } ) ^2 \} ^{-1/2} \right) ,
\\
G_p & = & e^{c\rho ^2 } \rho ^k \rho ^{\prime }  \{ \rho ^2 + (\rho ^{\prime } ) ^2 \} ^{-1/2} ,
\\
F^{\prime } & = & e^{c\rho ^2 } \rho ^{k+1} .
\end{eqnarray*}
By (\ref{rhocont}) and (\ref{cond11}), $G_p $ and $F^{\prime } $ are bounded on $(0, \pi)$. 
Moreover, since 
$P_{\mu } (\rho ) < +\infty $, we also have $G_r z \in L^1 ((0, \pi ))$. 
Integrating (\ref{Euler}) between $0$ and $t\in (0, \pi /2 )$ gives 
\begin{eqnarray}
\nonumber
-G_p (\rho (t), \rho ^{\prime } (t) )z(t) & = & \int _0 ^ t ( 
\gamma F^{\prime } - G_r )z \, d\theta 
\\
\nonumber
 & = & \int_0 ^t z e ^{c\rho ^2 } \rho ^k \left(  \gamma \rho - (2c\rho + (k/\rho )) 
\left[ \rho ^2 + (\rho ^{\prime } )^2 \right] ^{1/2} - \rho \left[ \rho ^2 + (\rho ^{\prime } )^2 \right] ^{-1/2} \right) \, d\theta 
\\
 & \leq & \gamma \int_0 ^{t } e^{c\rho ^2 } \rho ^{k+1 } z\, 
 d\theta \leq C \int _0 ^t z\, d\theta  \leq  Ct^{k+ 1 } .  
\label{intident}
\end{eqnarray}
On the other hand, if $\rho ^{\prime } (t) <0$, 
then (\ref{cond11}) and the boundedness of $\rho $ show that
\begin{eqnarray}
\label{estGp}
-G_p (\rho (t), \rho ^{\prime } (t) ) z(t)  & = & - e^{c(\rho (t) )^2  } 
(\rho (t))^k \rho^{\prime } (t) \left[  \rho ^2 (t) + 
(\rho ^{\prime } )^2 \right] ^{-1/2} \sin ^k t
\\
\nonumber
 & \geq & - C \rho ^{\prime } (t) \left[  C ^2  + (\rho ^{\prime } )^2 \right] ^{-1/2} t^k .
\end{eqnarray}
Furthermore, the estimate (\ref{cond6}) and the boundedness of $\rho $ imply that there is a constant $d_5 >0$ such that  
\begin{equation}
\label{rhoprimefromabove}
\rho ^{\prime } (t) \leq d_5 t \quad \forall t\in (0, \pi /2 ).
\end{equation}
Now (\ref{estGp}), (\ref{intident}) and (\ref{rhoprimefromabove}) imply that 
\begin{equation}
\label{rhoprimeovert}
\rho ^{\prime } (t)/t \ \mbox{ is bounded on $(0, \pi /2 )$.}
\end{equation}
In particular we have $\rho \in C^1 ([0, \pi ])$ and $\rho ^{\prime } (0) = \rho ^{\prime } (\pi )=0$. 

Finally, using (\ref{Euler}), a short calculation gives
\begin{equation}
\label{ident}
\gamma \rho   =  \frac{- \rho ^2 \rho ^{\prime \prime }  + 
\rho (\rho ^{\prime } )^2  }{( \rho ^2 + (\rho ^{\prime } )^2 )^{3/2} } 
 - \frac{ [(k/\rho ) + 2c\rho ] 
(\rho ^{\prime } )^2  + k 
\rho ^{\prime } \cot t }{( \rho ^2 + (\rho ^{\prime } )^2 )^{1/2} }  .
\end{equation}
By (\ref{cond11}), (\ref{rhocont}) and (\ref{rhoprimeovert}) this implies that
\begin{equation}
\label{rhotwoprimebdd}
\rho   ^{\prime \prime } \in L^{\infty } ((0, \pi )) .
\end{equation}
{\sl Step 6: $M$ is a half-disk }
\\[0.1cm]
Note first
that the derivatives $G_{rr} $, $G_{rp} $, $G_{pp} $ and $F^{\prime \prime } $ are bounded, in view of the properties (\ref{cond11}), (\ref{rhocont}) and (\ref{rhoprimeovert}).
 
Since $\rho $ is a minimizer of (\ref{isop2}), the second variation of $
P_{\mu }$ at $\rho $ in $K$ is nonnegative. This means that
\begin{equation}
0
\leq 
\int_0 ^{\pi }
\left( 
G_{rr}
\kappa ^2 
+2G_{rp}
\kappa \kappa ^{\prime }
+G_{pp} (\kappa ^{\prime \prime })^2
-\gamma F^{ \prime \prime } \kappa ^2 
\right)
z\,
d\theta ,  
\label{inequ22}
\end{equation}
for every $\kappa \in W^{1,2} ((0,\pi ))$ such that
\begin{equation}
\int_{0}^{\pi }F^{\prime }\kappa z\, d\theta =0.  
\label{F'kz=0}
\end{equation}
Furthermore, dividing (\ref{Euler}) by $z$ and then differentiating
yields
\begin{equation}
G_{rr}\rho ^{\prime }+G_{rp}\rho ^{\prime \prime }
- \frac{d}{d\theta }
\left( G_{rp}\rho ^{\prime }+G_{pp}\rho ^{\prime \prime }\right) -\left(
G_{pr}\rho ^{\prime }+
G_{pp} \rho ^{\prime \prime }\right) 
\frac{z^{\prime }}{z}
-G_p \left( \frac{z^{\prime }}{z}\right) ^{\prime }=\gamma F^{\prime
\prime }\rho ^{\prime }\quad \mbox{in }\ (0,\pi ).
\label{id3}
\end{equation}
Multiplying (\ref{id3}) by $\rho ^{\prime }z$ and then integrating by parts,
we obtain
\begin{equation}
\int_{0}^{\pi }G_{p}\rho ^{\prime }\left( \frac{z^{\prime }}{z}\right)
^{\prime } z\, d\theta =\int_0 ^{\pi }\left( G_{rr} (\rho ^{\prime
} )^2 +2G_{rp}\rho ^{\prime }\rho ^{\prime \prime }+G_{pp} (\rho ^{\prime
\prime })^2 -\gamma F^{\prime \prime } (\rho ^{\prime } )^2
\right) z \, d\theta .  \label{id4}
\end{equation}
Note that we may use (\ref{inequ22}) with 
 $\kappa =\rho ^{\prime }\in W^{1,\infty } ((0, \pi )) $. 
 This shows that the right-hand side of (\ref{id4}) is nonnegative. 
 On the other hand, 
\begin{equation}
\int_{0}^{\pi }G_p \rho ^{\prime }\left( \frac{z^{\prime }}{z}\right)
^{\prime }z\,d\theta
=
-k \int_{0}^{\pi } e ^{c \rho ^2 }
\rho ^k \left( \rho ^2 +(\rho ^{ \prime }) ^2 \right) ^{-1/2} 
(\rho ^{\prime }) ^2 \sin ^{k-2}\theta \, d\theta  \leq 0. 
\label{inequ3}
\end{equation}
Hence 
\begin{equation}
\label{id5}
\int_0 ^{\pi } e ^{ c\rho ^2 }
\rho ^k \left( \rho ^2 +(\rho ^{\prime } ) ^2 \right) ^{-1/2}
(\rho ^{\prime }) ^2 \sin ^{k-2}\theta \,d\theta  = 0,
\end{equation}
which implies that $\rho ^{\prime }=0$ in $[0,\pi ]$. This means that $\rho $ is
constant in $[0,\pi ]$, and the result follows. 
$\mbox{ } \hfill \Box $
\subsection{The N-dimensional case. Proof of Theorem 1.1}
We proceed by induction over the dimension $N$. 
Note that the result for $N=2 $ is Theorem \ref{Theorem 3.3. }
\\
Assume that the assertion holds true for sets in $\mathbb{R}^N $,
for some $N\geq 2 $,  and more precisely, for all measures of the
type
$$
d\mu = x_N ^k \mbox{ exp } \{ c |x|^2 \} \, dx
$$
where $k \geq 0$ and $c\geq 0$.
\\
We write $y= (x', x_N , x_{N+1 } ) $ for points in $\mathbb{R}^{N+1}
$, where $x' \in \mathbb{R}^{N-1} $, and $x_N , x_{N+1 } \in
\mathbb{R} $. Let a measure $\nu $ on
$$
\mathbb{R}^{N+1} _+ := \{ y= (x', x_N , x_{N+1} ) \in \mathbb{R}
^{N+1} :\, x_{N+1} >0\}
$$
be given by
$$
d\nu = x_{N+1} ^k \mbox{ exp } \{ c (|x'|^2 + x_N ^2 + x_{N+1} ^2 )  \} \, dy .
$$
We define two measures $\nu _1 $ and $\nu _2 $ by
\begin{eqnarray*}
d \nu _1 & = & \mbox{ exp } \{ c |x'|^2   \} \, dx', \quad \mbox{
and }\\
d \nu _2 & = & x_{N+1} ^k \mbox{ exp } \{ c ( x_N ^2 + x_{N+1} ^2 )
\} \, dx_N dx_{N+1 } ,
\end{eqnarray*}
and note that $d\nu = d\nu _1 d\nu _2 $.
\\
Let $M $ be a subset of $\mathbb{R}_{+}^{N+1}$ having finite and
positive $\nu $-measure.
\\
We define $2$-dimensional slices
$$
M(x') := \{ (x_N , x_{N+1 } ):\, (x', x_N , x_{N+1} ) \in M\} ,
\quad (x' \in \mathbb{R}^{N-1}  ).
$$
Let $ M' := \{ x' \in \mathbb{R}^{N-1} :\, 0< \nu_2 (M(x')) \} $,
and note that $\nu _2 (M(x')) < +\infty $ for a.e. $x' \in M'$. For
all those $x'$, let $H(x') $ be the half disc in
$\mathbb{R} ^2 _+ $ centered at $(0,0)$ with $\nu_2 (M(x')) = \nu _2
(H(x'))$. (For convenience, we put $H(x')=
\emptyset $ for all $x' \in M'$ with $\nu_2 (M(x')) =+\infty $.) By
Theorem \ref{Theorem 3.3. } we have
\begin{equation}
\label{isopnu2} P_{\nu_2 } (H(x')) \leq P_{\nu _2 }
(M(x')) \quad \mbox{for a.e. } \ x' \in M'.
\end{equation}
Let
$$
H:= \{ y=(x', x_N , x_{N+1} ) :\, (x_N , x_{N+1} )\in
H(x') ,\, x' \in M'\} .
$$
The product structure of the measure $\nu $ tells us that
\\
{\bf (i)} $ \nu (M)=\nu (H)$, and 
\\
{\bf (ii)} the isoperimetric property for slices, (\ref{isopnu2}),
carries over to $M$, that is,
\begin{equation}
\label{isopnu} P_{\nu } (H) \leq P_{\nu } (M),
\end{equation}
(see for instance Theorem 4.2 of \cite{BBMP}).
\\
Note again, the slice $ H (x') = \{ ( x_N , x_{N+1 } ) :\,
(x' , x_N , x_{N+1 } ) \in H \} $ is a half disc  $\{
(r\cos \theta , r \sin \theta ):\, 0 < r < R(x') ,\, \theta \in (0,
\pi ) \} $, with $0<R(x') <+\infty $, ( $x'\in M'$). Set
$$
K:= \{ (x' , r) : \, 0< r <
 R(x' ), \, x' \in M' \} , 
$$
and introduce a measure $\alpha $ on $\mathbb{R}_+ ^N $ by
$$
d\alpha := a_k r^{k+1 } \mbox{ exp } \{ c(|x'|^2 + r^2 )\} \, dx'
dr,
$$
where
$$
a_k := \int_0 ^{\pi } \sin ^k \theta \, d\theta = B\left(
\frac{k+1}{2}, \frac{1}{2} \right).
$$
An elementary calculation then shows that
\begin{equation*}
\nu (H)  =  \alpha (K),
\end{equation*}
and
\begin{equation*}
P_{\nu } (H)  =  P_{\alpha } (K).
\end{equation*}
Let $B_R $ denote the open ball in $\mathbb{R}^N $ centered at the
origin,
 with radius $R$, and choose $R>0 $ such that
$$ 
\alpha  (B_R \cap \mathbb{R} ^N _+ ) = \alpha ( K).
$$
By the induction assumption it follows that
\begin{equation}
\label{isopalpha} P_{\alpha } (B_R \cap \mathbb{R} ^N _+ ) \leq
P_{\alpha } (K).
\end{equation}
Finally, let $ M^{\bigstar } $ be the half ball in $\mathbb{R} ^{N+1 } _+ $
centered at the origin, with radius $R$,
$$ 
M^{\bigstar }   := \{ y= (x', x_N ,
x_{N+1 } ):\, |x'|^2 + x_N ^2 + x_{N+1 } ^2 < R^ 2 , \, x_{N+1 }
>0 \} .
$$
Then
\begin{equation*}
\nu (M^{\bigstar })  =  \nu (M)\quad \\
\end{equation*}
 and
\begin{equation*}
 P_{\nu }
(M^{\bigstar}) =  P_{\alpha } (B_R \cap \mathbb{R}^N _+ ).
\end{equation*}
Together with (\ref{isopalpha}) and (\ref{isopnu}) we
find
\begin{equation*}
\label{isopN+1} P_{\nu } (M^{\bigstar }) \leq P_{\nu } (M),
\end{equation*}
that is, the isoperimetric property holds for $N+1 $ in place of $N$
dimensions. The Theorem is proved. $\hfill \Box $
\section{Application to a class of degenerate elliptic equations}
\subsection{Notation and preliminary results}
First we introduce the notion of weighted rearrangement. 
For an exhaustive treatment of rearrangements we refer to
\cite{Ba}, \cite{CR}, \cite{H}, \cite{Ka}, \cite{Ke}. 
\\
Let the measure $\mu $ be given by (\ref{intro_dmu}), and let $M$ be
a measurable subset of $\mathbb{R}^N_+$. The distribution function
of a Lebesgue measurable function $u: M \rightarrow \mathbb{R} $,
with respect to $d\mu $, is the function $m_{\mu }$ defined by
\begin{equation*}
m_{\mu }(t)=\mu \left( \left\{ x\in M :\left\vert u(x)\right\vert
>t\right\} \right) ,\text{ }\forall t\geq 0.
\end{equation*}
The decreasing rearrangement of $u$ is the function $u^{\ast }$
 defined by
\begin{equation*}
u^{\ast }(s)=\inf \left\{ t\geq 0:
m_{\mu }(t)\leq s\right\} ,\text{ }\forall
s\in \left ( 0,\mu \left( M \right) \right] .
\end{equation*}
Let $C_{\mu }$ be the $\mu $-measure of $B_1\cap \mathbb{R}_{+}^{N}$,
that is,
\begin{equation*}
C_{\mu } = \frac{1}{2} (N-1) \omega _{_{N-1 }} B\left( \frac{k+1}{2} ,
\frac{N-1}{2} \right) ,
\end{equation*}
and let a function $ \psi
\left( r\right) $ be defined by
\begin{equation*}
\psi \left( r\right) =\int_{0}^{r}\exp \left( ct^{2}\right) t^{N+k-1}dt.
\end{equation*}
Let  $M ^{\bigstar }$ be defined as in Theorem 1.1, that is, 
\begin{equation}
M ^{\bigstar }=B_{r^{\bigstar }} \cap \mathbb{R}_{+}^{N},
\label{Omega-star-dv}
\end{equation}
where
\begin{equation}
r^{\bigstar }=\psi ^{-1}\left(  \frac{ \mu (M )}{C_{\mu } } \right) .
\label{R-star-dv}
\end{equation}
The rearrangement $u^{\bigstar }$ of $u$, by its definition given in \eqref{defrearr},  is 
\begin{equation*}
u^{\bigstar }(x)=u^{\ast }\left( C_{\mu }\psi \left( \left\vert x\right\vert
\right) \right) ,\text{ }\forall x\in M ^{\bigstar }.
\end{equation*}
The isoperimetric inequality in Theorem \ref{Theorem 3.2. }
can be also stated as follows
\begin{equation*}
P_{\mu }(M )\geq I_{\mu }(\mu (M )),
\end{equation*}
where $I_{\mu }(\tau )$ is the function such that $P_{\mu }(M
^{\bigstar })=I_{\mu }(\mu (M ^{\bigstar })),$ or equivalently
\begin{equation}
I_{\mu }(\tau )=C_{\mu }\exp \left( c\left[ \psi ^{-1}\left( \frac{\tau
}{C_{\mu }} \right) \right] ^{2}\right) \left[ \psi ^{-1}\left(
\frac{\tau }{C_{\mu }} \right) \right] ^{N+k-1}.
\label{Def-Isop-Funct-dv}
\end{equation}
The fact that half balls $B_R \cap \mathbb{R}^N _+ $ are 
isoperimetric for the weighted measure $\mu $ imply a
Polya-Szeg\"{o} - type inequality (see
\cite{T}, p. 125).  

\begin{theorem}\label{Theorem 4.1. } 
{\sl Let $D$ be an open set with finite $\mu $--measure, and   
let the space $V^2(D ,d\mu )$ be given by Definition \ref{Definition 2.5. } Then we have for every 
function $u\in V^2 (D, d\mu )$, 
\begin{equation}
\int_{D }|\nabla u|^{2}d\mu \geq \int_{D ^{\bigstar }}|\nabla u^{\bigstar
}|^{2}d\mu .  
\label{Polya_Szego}
\end{equation}
}
\end{theorem}

Since rearrangements preserves the $L^{p\text{ }}$
norms, we have that the Rayleigh-Ritz quotient decreases under rearrangement
i.e.
\begin{equation*}
\frac{\int_{D }|\nabla u|^{2}d\mu }{\int_{D }u^{2}d\mu }\geq \frac{
\int_{D ^{\bigstar }}| \nabla u^{\bigstar }|^{2}d\mu }{\int_{D
^{\bigstar }}\left( u^{\bigstar }\right) ^{2}d\mu },\text{ }\
\forall u\in V^2(D ,d\mu ).
\end{equation*}
The following Poincar\`{e} type inequality states the continuous
embedding of $V^2(D ,d\mu )$ in $L^{2}(D ,d\mu )$. It is a consequence of 
some one-dimensional inequalities (see
\cite{Ma}, Theorem 2, p. 40).

\begin{corollary}\label{Corollary 4.2. } {\sl Let $D$ be an open subset of $\mathbb{R} ^N _+ $. 
Then there exists a constant $C$, such that for every $u\in V^2(D ,d\mu )$,
\begin{equation*}
\int_{D }u^{2}d\mu \leq C\int_{D }\left\vert \nabla u\right\vert
^{2}d\mu .
\end{equation*}
}
\end{corollary}

\subsection{Comparison result}
Now we are in a position to obtain sharp estimates for the solution
to problem (\ref{P}). By a weak solution to such a problem we mean a
function $u$ belonging to $V^2(D ,d\mu )$ such that
\begin{equation}
\int_{D }A(x)\nabla u\nabla \chi d\mu =\int_{D }f \chi d\mu ,
\label{wsol-dv}
\end{equation}
for every $\chi \in C^{1}(\bar{D})$ such that $\chi =0$ on the set
$\partial D \setminus \left\{ x_{N}=0\right\} $.
\\[0.2cm]
\noindent {\sl Proof of Theorem 1.2 : }  Note first that the existence
of a unique solution to problems (\ref{P}) and (\ref{Symm_Probl_dv})
is ensured by the Lax - Milgram Theorem. Arguing as in \cite{T2} 
(see for instance \cite{BCM}, p. 363), 
we get
\begin{equation}
1\leq \left\{ \left[ I_{\mu }\left( m_{u}(t)\right) \right]
^{-2}\int_{0}^{m_{u}(t)}f^{\ast }(\sigma )d\sigma \right\} \left(
-m_{u}^{\prime }(t)\right)  \label{e-dv}
\end{equation}
and
\begin{equation}
u^{\ast }\left( s\right) \leq \int_{s}^{\mu (D )}\left( I_{\mu
}^{-2}(l)\int_{0}^{l}f^{\ast }(\sigma )d\sigma \right) dl,  \label{B-dv}
\end{equation}
Using (\ref{Def-Isop-Funct-dv}) in (\ref{B-dv}) we obtain
\begin{align*}
u^{\bigstar }\left( x\right) & \leq \frac{1}{{C_{\mu }}^{2}}\int_{C_{\mu }\psi
\left( \left\vert x\right\vert \right) }^{\mu (D )}\left\{ \exp \left(
-2c\left[ \psi ^{-1}\left( \frac{l}{C_{\mu }}\right) \right] ^{2}\right) \left(
\psi ^{-1}\left( \frac{l}{C_{\mu }}\right) \right)
^{-2N-2k+2}\int_{0}^{l}f^{\ast }(\sigma )d\sigma \right\} dl \\
& =\frac{1}{{C_{\mu }}}\int_{\left\vert x\right\vert }^{r^{\bigstar }}\exp
\left( -c\eta ^{2}\right) \eta ^{-N-k+1}\left( \int_{0}^{C_{\mu }\psi \left(
\eta \right) }f^{\ast }(\sigma )d\sigma \right) d\eta \text{ \ \ \ \ \ \ \ \
\ \ \ \ \ \ \ \ \ }\left( \eta :=\psi ^{-1}\left( \frac{l}{C_{\mu }}\right)
\right) \\
& =\int_{\left\vert x\right\vert }^{r^{\bigstar }}\exp \left( -c\eta
^{2}\right) \eta ^{-N-k+1}\left( \int_{0}^{\eta }f^{\ast }(C_{\mu }\psi (\xi
))\xi ^{N+k-1}\exp \left( c\xi ^{2}\right) d\xi \right) d\eta \text{ \ \ \ \
}\left( \sigma :=C_{\mu }\psi (\xi )\right) \\
& =\int_{\left\vert x\right\vert }^{r^{\bigstar }}\exp \left( -c\eta
^{2}\right) \eta ^{-N-k+1}\left( \int_{0}^{\eta }f^{\bigstar }(\xi )\xi
^{N+k-1}\exp \left( c\xi ^{2}\right) d\xi \right) d\eta \\
& =w(x).
\end{align*}
Now let us show (\ref{Grad-Est-dv}). Arguing as in \cite{BCM}, p.
363--364 (see also \cite{T2}), we derive
\begin{eqnarray*}
-\frac{d}{dt}\int_{\left\vert u\right\vert >t}\left\vert \nabla u\right\vert
^{q}d\mu &\leq &\left( \int_{0}^{m_{u}(t)}f^{\ast }(s)ds\right) ^{q/2}\left(
-m_{u}^{\prime }(t)\right) ^{1-q/2} \\
&\leq &\left( I(m_{u}(t))\right) ^{-q}\left( \int_{0}^{m_{u}(t)}f^{\ast
}(s)ds\right) ^{q}\left( -m_{u}^{\prime }(t)\right) .
\end{eqnarray*}
Integrating the last inequality between $0$ and $+\infty $, we get
\begin{eqnarray*}
\int_{D }\left\vert \nabla u\right\vert ^{q}d\mu &=&\int_{0}^{+\infty }\left[
-\frac{d}{dt}\int_{\left\vert u\right\vert >t}\left\vert \nabla u\right\vert
^{q}d\mu \right] dt \\
&\leq &\int_{0}^{+\infty }\left( I_{\mu }(m_{u}(t))\right) ^{-q}\left(
\int_{0}^{m_{u}(t)}f^{\ast }(\sigma )d\sigma \right) ^{q}\left(
-m_{u}^{\prime }(t)\right) dt \\
&\leq &\int_{0}^{\mu (D )}\left( I_{\mu }(s)\right) ^{-q}\left(
\int_{0}^{s}f^{\ast }(\sigma )d\sigma \right) ^{q}ds.
\end{eqnarray*}
Now a straightforward calculation yields
\begin{eqnarray*}
 & & \int_{D }\left\vert \nabla u\right\vert ^{q}d\mu 
\leq 
C_{\mu }\int_{0}^{\mu (D )}\exp \left( -qc\left[ \psi ^{-1}\left(
\frac{s}{C_{\mu }}\right)
\right] ^{2}\right) \left[ \psi ^{-1}\left(
\frac{s}{C_{\mu}}\right) \right]
^{-q(N+k-1)}\left( \int_{0}^{s}f^{\ast }(\sigma )d\sigma \right) ^{q}ds 
\\
 & = & 
C_{\mu }^{2}\int_{0}^{R^{\star }}\exp \left( -qc\eta ^{2}\right) \eta
^{-q(N+k-1)}\left( \int_{0}^{C_{\mu }\psi (\eta )}f^{\ast }(\sigma )d\sigma
\right) ^{q}\exp \left( c\eta ^{2}\right) \eta ^{N+k-1}d\eta \\
 & = & 
C_{\mu }^{2}\int_{0}^{R^{\star }}\exp \left( (1-q)c\eta ^{2}\right)
\eta ^{(1-q)(N+k-1)}\left( \int_{0}^{C_{\mu }\psi (\eta )}f^{\ast
}(\sigma )d\sigma \right) ^{q}d\eta
\\
 & = & 
C_{\mu }^{2}\int_{0}^{R^{\star }}\left( \int_{0}^{\eta }f^{\ast
}(C_{\mu }\psi (\rho ))C_{\mu }\exp \left( c\rho ^{2}\right) \rho
^{N+k-1}d\rho \right)
^{q}\exp \left( (1-q)c\eta ^{2}\right) \eta ^{(1-q)(N+k-1)}d\eta \\
 & = & 
 C_{\mu }^{2+q}\int_{0}^{R^{\star }}\left( \int_{0}^{\eta }f^{\bigstar }(\rho
)\exp \left( c\rho ^{2}\right) \rho ^{N+k-1}d\rho \right) ^{q}\exp \left(
(1-q)c\eta ^{2}\right) \eta ^{(1-q)(N+k-1)}d\eta \\
 & = &
\int_{D }\left\vert \nabla w\right\vert ^{q}d\mu .
\end{eqnarray*}
{\bf Acknowledgement: }
The first author wants to thank the University of Naples for a visiting appointment.

\end{document}